\documentclass[11pt]{article}
\usepackage{tikz}
\usepackage[all]{xy}
\usepackage{amsfonts}
\usepackage{amsfonts}
\usepackage{amsfonts}
\usepackage{mathrsfs}
\usepackage{shuffle}
\usepackage{bm}
\usepackage{cite}
\usepackage[colorlinks,linkcolor=blue,citecolor=blue]{hyperref}
\usepackage{amssymb,amsmath} 
\usepackage{makecell}
\usepackage{tikz}
\usepackage[caption=false]{subfig}
\usepackage[utf8]{inputenc}
 \allowdisplaybreaks
\usetikzlibrary{arrows,shapes,chains}
 \allowdisplaybreaks

\def\kk{k_1,\ldots,k_r}
\catcode`!=11
\let\!int\int \def\int{\displaystyle\!int}
\let\!lim\lim \def\lim{\displaystyle\!lim}
\let\!sum\sum \def\sum{\displaystyle\!sum}
\let\!sup\sup \def\sup{\displaystyle\!sup}
\let\!inf\inf \def\inf{\displaystyle\!inf}
\let\!cap\cap \def\cap{\displaystyle\!cap}
\let\!max\max \def\max{\displaystyle\!max}
\let\!min\min \def\min{\displaystyle\!min}
\let\!frac\frac \def\frac{\displaystyle\!frac}
\catcode`!=12

\let\oldsection\section
\renewcommand\section{\setcounter{equation}{0}\oldsection}

\allowdisplaybreaks
\def\pf{\it{Proof.}\rm\quad}

\def\N{\mathbb{N}}
\def\A{{\rm A}}

\def\su{\sum\limits_{n=1}^\infty}

\newcommand\bfk{{\bf k}}
\newcommand\bfl{{\bf l}}

\newtheorem{thm}{Theorem}[section]
\newtheorem{defn}[thm]{Definition}
\newtheorem{lem}[thm]{Lemma}
\newtheorem{cor}[thm]{Corollary}

\newtheorem{re}[thm]{Remark}

\newtheorem{pro}[thm]{Proposition}
\begin{document}
\title {\bf Some results on Arakawa-Kaneko, Kaneko-Tsumura functions and related functions}
\author{
{Maneka Pallewatta$^{a,}$\thanks{Email: maneka.osh@gmail.com} \quad Ce Xu$^{a,b,}$\thanks{Email: cexu2020@ahnu.edu.cn}}\\[1mm]
\small a. Graduate School of Mathematics, Kyushu University, Motooka\\
\small  Nishi-ku, Fukuoka 819-0389, Japan\\
\small b.  School of Mathematics and Statistics, Anhui Normal University,\\ \small Wuhu 241002, P.R. China\\
[5mm]
Dedicated to Professor Masanobu Kaneko on the occasion of his 60th birthday}

\date{}
\maketitle \noindent{\bf Abstract} Recently, the level two analogue of multiple polylogarithm function ${\rm A}(k_1,\ldots,k_r;z)$ and Arakawa-Kaneko zeta function $\psi(k_1,\ldots,k_r;s)$ were introduced by M.~Kaneko and H.~Tsumura, for $\kk \in \mathbb{Z}_{\ge 1}$ . In this paper, we investigate some of their  special relations. In particular, we prove some explicit forms of ${\rm A}(k_1,\ldots,k_r;z)$ and $\psi(k_1,\ldots,k_r;s)$. Also, we introduce a level $m$ anlogue of the Arakawa-Kaneko zeta functions.
\\[2mm]
\\[2mm]

\section{Introduction}
We begin with some basic notations. Let us consider the positive index set $ {\bfk_r}:= (k_1,\ldots, k_r)$. The quantities
$|{\bfk_r}|:=k_1+\cdots+k_r$ and $ {\rm dep}({\bfk_r}):=r,$
are called the weight and depth of ${\bfk_r}$, respectively. If $k_r>1$, ${\bfk_r}$ is called \emph{admissible}. For ${\bfk_r}:= (k_1,\ldots, k_r)$, set ${\bfk_0}:=\emptyset$, $({\bfk_r})_{+}:=(k_1,\ldots,k_{r-1},k_r+1)$ and $({\bfk_r})_{-}:=(k_1,\ldots,k_{r-1},k_r-1)$.

The subject of this paper is the level two analogue of Arakawa-Kaneko and related functions, which is a generalisation of the single-variable multiple zeta function.

The \emph{Arakawa-Kaneko zeta function} (\cite{AM1999}) is defined by
\begin{align}\label{a1}
\xi(k_1,k_2\ldots,k_r;s):=\frac{1}{\Gamma(s)} \int\limits_{0}^\infty \frac{t^{s-1}}{e^t-1}{\rm Li}_{k_1,k_2,\ldots,k_r}(1-e^{-t})dt,
\end{align}
for $k_1,k_2,\ldots,k_{r} \in \mathbb{Z}_{\ge 1}$ and $\Re(s)>0$, where ${\mathrm{Li}}_{{{k_1},{k_2}, \cdots ,{k_r}}}\left( z \right)$ is the \emph{multiple poly-logarithm} defined by
\begin{align}\label{a2}
&{\mathrm{Li}}_{{{k_1},{k_2}, \cdots ,{k_r}}}\left( z \right): = \sum\limits_{1 \le {n_1} <  \cdots  < {n_r}} {\frac{{{z^{{n_r}}}}}{{n_1^{{k_1}}n_2^{{k_2}} \cdots n_r^{{k_r}}}}}\quad z \in \left[ { - 1,1} \right).
\end{align}

In a recent paper \cite{KT2018}, Kaneko and Tsumura introduced and studied a new kind of Arakawa-Kaneko type functions
\begin{align}\label{Def-eta}
\eta(k_1,k_2\ldots,k_r;s):=\frac{1}{\Gamma(s)} \int\limits_{0}^\infty \frac{t^{s-1}}{1-e^t}{\rm Li}_{k_1,k_2,\ldots,k_r}(1-e^{t})dt.
\end{align}
We call them \emph{Kaneko-Tsumura $\eta$-function}.

In the past two decades, the study of Arakawa-Kaneko and related functions has attracted the attention of many mathematicians. Apart from the actual evaluation of the functions, one of the main questions that one sets out to solve is that whether or not Arakawa-Kaneko zeta functions can be expressed in terms of a linear combination of \emph{multiple zeta values} (MZVs) (\cite{H1992,DZ1994})
\begin{align}\label{a4}
\zeta(k_1,\ldots,k_{r-1},k_r):=\sum\limits_{0<n_1<\cdots<n_r}\frac{1}{n_1^{k_1}\cdots n_{r-1}^{k_{r-1}} n_r^{k_r}},
\end{align}
and \emph{single-variable multiple zeta functions}
\begin{align}\label{a6}
\zeta(k_1,\ldots,k_{r-1},s):=\sum\limits_{0<n_1<\cdots<n_r}\frac{1}{n_1^{k_1}\cdots n_{r-1}^{k_{r-1}} n_r^{s}},
\end{align}
where $\bfk_r$ is an admissible index and $\Re(s)>1$.

In \cite{AM1999,KT2018,Ku2010}, for any index ${\bf k}_r:=(k_1,\ldots,k_r)$, the special values $\xi({\bf k}_r;s)$ at positive integers are analytically computed and written in terms of multiple zeta values.
Further, M.~Kaneko and H.~Tsumura \cite{KT2018} also proved that for a general index set ${\bf k}_r$, the function $\xi({\bf k}_r;s)$ can also be expressed by multiple zeta functions (but not explicit formulas). For example,
\begin{align}\label{a7}
\xi(\{1\}_{r-1},k;s)=&(-1)^{k-1}\sum\limits_{a_+\cdots+a_k=r \atop \forall a_j\geq 0} \binom{s+a_k-1}{a_k} \zeta(a_1+1,\ldots,a_{k-1}+1,a_k+s)\nonumber\\
&+\sum\limits_{j=0}^{k-2}(-1)^j \zeta(\{1\}_{r-1},k-j)\zeta(\{1\}_{j},s).
\end{align}
Some related results (e.g. duality formulas etc.) for Arakawa-Kaneko type functions can be found in the works of \cite{KTB2018,KO2018,Y2016}.
Here, $\{l\}_m$ denotes the sequence $\underbrace{l,\ldots,l}_{m \text{\;-times}}$.

Recently, Kaneko and Tsumura defined the single variable multiple zeta function of level two as follows.
\begin{defn}\label{def:1}(Kaneko, Tsumura \cite{KTA2018})
For $k_1,\ldots,k_{r-1} \in \mathbb{Z}_{\ge 1}$ and $\Re{(s)}>1$, we write
\begin{equation}\label{a10}
    T_0(k_1,\ldots,k_{r-1},s)=\sum_{\substack{0<m_1<\cdots<m_r \\ m_i\equiv i \pmod{2}}} \frac{1}{m_1^{k_1}\cdots m_{r-1}^{k_{r-1}}m_r^s}.
\end{equation}
Furthermore, as its normalized version,
\begin{equation}\label{a11}
    T(k_1,\ldots,k_{r-1},s)=2^r  T_0(k_1,\ldots,k_{r-1},s),
\end{equation}
which is called multiple $T$-function.
 The values $T(k_1,\ldots,k_{r-1},k_r)$ ($k_j \in \mathbb{Z}_{\ge 1}$, $k_r\ge 2:$ admissible) are called the multiple T-values (MTVs).
\end{defn}

Using these functions, Kaneko and Tsumura defined a level two analogue of $\xi(k_1,\ldots,k_r;s)$ which we call Kaneko-Tsumura function as follows.
\begin{defn}(Kaneko, Tsumura \cite[\S 5]{KTA2018})\label{df-3.3}
For index $\bfk_r$ and $\Re{(s)}>0$, we write
\begin{equation}\label{a8}
    \psi(\kk;s)=\frac{1}{\Gamma(s)}\int_0^\infty t^{s-1} \frac{{\rm A}(\kk;\tanh t/2)}{\sinh(t)} dt,
\end{equation}
where
\begin{align}\label{a9}
&{\rm A}(k_1,k_2,\ldots,k_r;z): = 2^r\sum\limits_{1 \le {n_1} <  \cdots  < {n_r}\atop n_i\equiv i\ {\rm mod}\ 2} {\frac{{{z^{{n_r}}}}}{{n_1^{{k_1}}n_2^{{k_2}} \cdots n_r^{{k_r}}}}},\quad z \in \left[ { - 1,1} \right).
\end{align}
\end{defn}
In particular, if $s\in \mathbb{Z}_{\ge 1}$ then we call (\ref{a8}) the Kaneko-Tsumura $\psi$-values. Here, ${\rm A}(k_1,\ldots,k_r;z)$ is $2^r$ times ${\rm Ath}(k_1,\ldots,k_r;z)$ which was introduced in \cite[\S5]{KTA2018}.
When $k_r >1$, we see that
\begin{equation*}
    \A(\kk;1)=T(\kk).
\end{equation*}

${\rm A}(\mathbf{k}_r;z)$ satisfies the shuffle relation corresponding to multiple zeta values (see \cite{H1992} and \cite{HO2003}). Let $\mathfrak{H} := \mathbb{Q} \langle x,y \rangle$
be the non-commutative polynomial ring in two indeterminates $x$ and $y$. We refer to monomials in $x$ and $y$ as words. We also define subrings, $\mathfrak{H}^1 := \mathbb{Q} +y \mathfrak{H} $
and $\mathfrak{H}^0 := \mathbb{Q} + y \mathfrak{H} x$.
For any integer $k >0$, put $z_k=yx^{k-1}$. Then the ring $\mathfrak{H}^1$ is freely generated by $z_k$ $(k \ge 1)$. When $k \ge 2$, $z_k$ is contained in $\mathfrak{H}^0$. But $\mathfrak{H}^0$ is not freely generated by $z_k$ $(k \ge 2)$.
Now let us define the evaluation map $Z:\mathfrak{H}^0 \to \mathbb{R} $ by
\begin{equation}
    Z(z_{k_1} \cdots z_{k_r}) := {\rm A}(k_1, \ldots , k_r;z).
\end{equation}
We define the shuffle product $\shuffle$ on $\mathfrak{H}$ inductively by
\begin{align*}
     1\shuffle w &=w\shuffle1 =w \\
    u_1 w_1\shuffle u_2 w_2 &=u_1(w_1\shuffle u_2 w_2)+u_2(u_1 w_1 \shuffle w_2)
\end{align*}
for any words $w,w_1,w_2 \in \mathfrak{H}$ and $u_1,u_2 \in \{x,y\}$, with $\mathbb{Q}$-bilinearity.
The shuffle product is commutative and associative.

The main purpose of this paper is to study the functions $\psi({\bf k}_r;s), \A({\bf k}_r;z)$ and MTVs. In particular, we give some explicit formulas of $\psi({\bf k}_r;s)$ in terms of multiple $T$-functions.

In $\S2$, we first obtain new formula for $\psi({\bf k}_r;z)$ corresponding to $\xi({\bf k}_r;s)$ function. Secondly, we prove some explicit forms of new identities for $\A({\bf k}_r;z)$ and $\psi({\bf k}_r;s)$ by using the methods of iterated integral representations of series. Similarly, we can obtain formulas for ${\rm Li}_{{\bf k}_r}(z)$ and $\xi({\bf k}_r;s)$.

In $\S3$, we discuss a general duality relation of Kaneko-Tsumura $\psi$-values and give an explicit formula.

In $\S4$, we introduce and study a ``level $m$" analogue of the multiple polylogarithm and the multiple zeta functions. Moreover, we give an equation system and prove that. If the solution of this equation system exists, then the level $m$ analogue of $\xi(k_1,\ldots,k_r;s)$ can be defined. Furthermore, we can deduce many results corresponding to $\xi(k_1,\ldots,k_r;s)$ and $\psi(k_1,\ldots,k_r;s)$.

\section{Main results and proofs}

\subsection{Relations among the functions $\psi$ and $T$}
In this section, we present our main results on the Kaneko-Tsumura zeta functions. We deduce that the Kaneko-Tsumura zeta function can be written as a linear combination of multiple $T$-functions.

In order to prove the main results, we establish the following lemmas.

\begin{lem}\label{e1}(\cite{KT2018}, cf. \cite{AM1999}) {\rm (i)} For index $\bfk_r$,
\begin{align}\label{b3}
\frac{d}{dz}{\mathrm{A}}({{k_1}, \cdots ,k_{r-1},{k_r}}; z)= \left\{ {\begin{array}{*{20}{c}} \frac{1}{z} {\mathrm{A}}({{k_1}, \cdots ,{k_{r-1}},{k_r-1}};z)
   {\ \ (k_r\geq 2),}  \\
   {\frac{2}{1-z^2}{\mathrm{A}}({{k_1}, \cdots ,{k_{r-1}}};z)\;\;\;\ \ \ (k_r = 1).}  \\
\end{array} } \right.
\end{align}
{\rm (ii)} For  $r \ge 1$,
  \begin{align}\label{b4}
{\rm A}({\{1\}_r};z)=\frac{1}{r!}({\rm A}(1;z))^r=\frac{(-1)^r}{ r!}\log^r\left(\frac{1-z}{1+z}\right).
\end{align}
\end{lem}
Similarly, we can obtain the following lemma.
\begin{lem}\label{e2}
{\rm (i)} For index $\bfk_r$,
     \begin{align}
         \frac{d}{dz} {\rm A}\left(\kk;\frac{1-z}{1+z}\right)=\left \{
         \begin{array}{c}
         -\frac{2}{1-z^2} {\rm A}\left(k_1,\ldots,k_{r-1},k_r-1;\frac{1-z}{1+z}\right) \; \; \;\;\; \; \;\;\;  (k_r \ge 2) \\
         -\frac{1}{z} {\rm A}\left(k_1,\ldots,k_{r-1};\frac{1-z}{1+z}\right) \;\;\; \; \;\;\; \; \;\;\;\;\; \; \;\;\; \; \;\;\;\;\; \; \;\;\; \; (k_r = 1).
         \end{array}
         \right.
     \end{align}
{\rm (ii)} For $r \ge 1$ \begin{align}
        {\rm A}\left(\{1\}_r;\frac{1-z}{1+z}\right)=\frac{1}{r!}{\rm A}^r\left(1;\frac{1-z}{1+z}\right)=\frac{(-1)^r}{r!}\log^r z.&&
    \end{align}
\end{lem}

\begin{lem}\label{lem-4.6}
For any index ${\bf{k}}_r$, we have
\begin{align}\label{eq-4.4}
   \frac{2}{1-z^2}{\rm A}\left(\{1\}_j;\frac{1-z}{1+z}\right) {\rm A}({\bf{k}}_r;z)=\frac{d}{dz} \left(\sum_{i=0}^j {\rm A}\left(\{1\}_{j-i};\frac{1-z}{1+z}\right) {\rm A}({\bf{k}}_r,i+1;z)\right).
\end{align}
\end{lem}
\pf
By using Lemma \ref{e1} and Lemma \ref{e2}, we can easily obtain the desired result by induction on $j$. \hfill$\square$

Similar to the Euler-type connection formula of multiple polylogarithm functions in \cite[Lem. 3.5]{KT2018},  we obtain the following formula associated with the level two analogue, for ${\rm A}({\bf k}_r;z)$.
\begin{thm}\label{thm-Ath} Let ${\bf k}_r$ be any index. Then we have
\begin{equation*}
    {\rm A}\left(\mathbf{k}_r;\frac{1-z}{1+z}\right)=\sum_{\mathbf{k}',j \ge 0} C_{\mathbf{k}_r}(\mathbf{k}';j) {\rm A}\left(\{1\}_j;\frac{1-z}{1+z}\right) {\rm A}(\mathbf{k}';z),
\end{equation*}
where the sums on the right runs over indices ${\bf k'}$ and integers $j\geq 0$ that satisfy $|{\bf k'}|+j\leq |{\bf k}_r|$, and $C_{{\bf k}_r}({\bf k'};j)$ is a $\mathbb{Q}$-linear combination of multiple $T$-values of weight $|{\bf k}_r|- |{\bf k'}|-j$. We understand ${\rm A}_{\emptyset}(z)=1$ and $|\emptyset|=0$ for the empty index $\emptyset$, and the constant $1$ is regarded as a multiple $T$-value of weight $0$.
\end{thm}
\pf
We prove this by induction on the weight $\mathbf{k}_r$. When $\mathbf{k}_r=(1)$, the trivial identity
\begin{equation*}
    \A_1\left(\frac{1-z}{1+z}\right)=\A_1\left(\frac{1-z}{1+z}\right)
\end{equation*}
itself gives the desired form, thus $C_{(1)}(\emptyset;0)=C_{(1)}((1);0)=0$ and $C_{(1)}(\emptyset;1)=1$.
Suppose the weight $|\mathbf{k}_r|>1$ and assume the statement holds for any index of weight less than $|\mathbf{k}_r|$.

For $\mathbf{k}_r=(\kk)$, set $(\mathbf{k}_r)_-=(k_1,\ldots,k_{r-1},k_r-1)$.

First, assume that $\mathbf{k}_r$ is admissible. Then by the differential relation and the induction hypothesis, we get
\begin{align}\label{eq-2-2}
    \frac{d}{dz}{\rm A}\left(\mathbf{k}_r;\frac{1-z}{1+z}\right)&=-\frac{2}{1-z^2}{\rm A}\left((\mathbf{k}_r)_-;\frac{1-z}{1+z}\right) \nonumber \\
    &=-\frac{2}{1-z^2}\sum_{\mathbf{l},j \ge 0} C_{(\mathbf{k}_r)_-}(\mathbf{l};j) {\rm A}\left(\{1\}_j;\frac{1-z}{1+z}\right) {\rm A}(\mathbf{l};z).
\end{align}
Let the depth of $\mathbf{l}$ be $s$. By substituting (\ref{eq-4.4}) from Lemma \ref{lem-4.6} into (\ref{eq-2-2}) and integrating, we get
\begin{equation*}
    {\rm A}\left(\mathbf{k}_r;\frac{1-z}{1+z}\right)=-\sum_{\mathbf{l} ,j \ge 0} C_{(\mathbf{k}_r)_-}(\mathbf{l};j) \left(\sum_{i=0}^j {\rm A}\left(\{1\}_{j-i};\frac{1-z}{1+z}\right) {\rm A}(\mathbf{l},i+1;z)\right)+C,
\end{equation*}
where $C$ is a constant. Since
\begin{equation*}
   \lim_{z\to 0} {\rm A}\left(\{1\}_{j-i};\frac{1-z}{1+z}\right) {\rm A}(\mathbf{l},i+1;z)=0,
\end{equation*}
we have $C=T(\mathbf{k}_r)$. Now we can obtain the desired result.

In order to prove the non-admissible case, we recall that ${\rm A}\left(\mathbf{k}_r;\frac{1-z}{1+z}\right)$ satisfies the shuffle relation (cf. \cite{IKZ2006}). Suppose $\mathbf{k}_r$ is not admissible. Then, we can write ${\rm A}\left(\mathbf{k}_r;\frac{1-z}{1+z}\right)$ as a polynomial of ${\rm A}\left(1;\frac{1-z}{1+z}\right)$ with each coefficient of ${\rm A}^i\left(1;\frac{1-z}{1+z}\right)$ being a linear combination of ${\rm A}\left({\mathbf{k}}';\frac{1-z}{1+z}\right), \mathbf{k}':\rm{admissible}$. Write this polynomial as
\begin{equation*}
    {\rm A}\left(\mathbf{k}_r;\frac{1-z}{1+z}\right)=\sum_{j=0}^m a_i\cdot{\rm A}^j\left(1;\frac{1-z}{1+z}\right) .
\end{equation*}
Then $a_i$ can be written in the desired form (admissible case). We know that
\begin{equation*}
   {\rm A}^j\left(1;\frac{1-z}{1+z}\right)=j! {\rm A}\left(\{1\}_j;\frac{1-z}{1+z}\right)
\end{equation*}
and
\begin{equation*}
    {\rm A}\left(\{1\}_i;\frac{1-z}{1+z}\right) {\rm A}\left(\{1\}_j;\frac{1-z}{1+z}\right) =\binom{i+j}{i}{\rm A}\left(\{1\}_{i+j};\frac{1-z}{1+z}\right).
\end{equation*}
Hence $a_i\cdot{\rm A}^j\left(1;\frac{1-z}{1+z}\right)$ can be written in the claimed form, and the proof is done. \hfill$\square$

We obtain a level two version of \cite[Proposition 2]{AM1999} which will be needed in proving our main results as follows.
\begin{pro}\label{prop-2-1}
\begin{enumerate}
     \item For $\Re {(s)}>1$
     \begin{equation*}
    T(k_1,\ldots,k_{n-1},s)=\frac{1}{\Gamma(s)}\int_0^\infty \frac{t^{s-1}}{\sinh(t)}{\A(k_1,\ldots,k_{n-1};e^{-t})} dt.
\end{equation*}

\item For $\Re {(s)}>1, n\ge2, j\ge0$
  \begin{equation*}
    \int_0^\infty t^{s+j-1}{\A(k_1,\ldots,k_{n-1};e^{-t})} dt=\Gamma(s+j)T(k_1,\ldots,k_{n-2},s+j+k_{n-1}).
\end{equation*}
\end{enumerate}
\end{pro}
The proof is similar to the proof of \cite[Proposition 2]{AM1999}. Therefore, we omit the proof.

From Theorem \ref{thm-Ath}, we can obtain formulas expressing $\psi(\mathbf{k}_r;s)$ in terms of multiple $T$-zeta functions.
\begin{thm}\label{thm-psi}
Let $\mathbf{k}_r$ be any index set. The function $\psi(\mathbf{k}_r;s)$ can be written in terms of multiple $T-$functions as
\begin{equation*}
  \psi(\mathbf{k}_r;s)=\sum_{\mathbf{k}',j \ge 0} C_{\mathbf{k}_r}(\mathbf{k}';j)\binom{s+j-1}{j}T(\mathbf{k}';s+j)
\end{equation*}
Here, the sum is over indices $\mathbf{k}'$ and integers $j\ge0$ that satisfy $|\mathbf{k}'|+j\le |\mathbf{k}_r|$, and $C_{\mathbf{k}_r}(\mathbf{k}';j)$ is the same as in Theorem \ref{thm-Ath}.
\end{thm}
\pf
Let $r, \; l$ be the depths of $\mathbf{k}_r$ and $\mathbf{k}'$ respectively.
Put $z=e^{-t}$ in Theorem \ref{thm-Ath}.

\begin{equation*}
    {\rm A}\left(\mathbf{k}_r;\frac{1-e^{-t}}{1+e^{-t}}\right)=\sum_{\mathbf{k}',j \ge 0} C_{\mathbf{k}_r}(\mathbf{k}';j) {\rm A}\left(\{1\}_j;\frac{1-e^{-t}}{1+e^{-t}}\right) {\rm A}(\mathbf{k}';e^{-t}).
\end{equation*}
By using Lemma \ref{e1} we can write the above equation as
\begin{equation}\label{eq-2-3}
  {\rm A}\left(\mathbf{k}_r;\tanh t/2\right)=\sum_{\mathbf{k}',j \ge 0} C_{\mathbf{k}_r}(\mathbf{k}';j) \frac{t^j}{j!} {\rm A}(\mathbf{k}';e^{-t}).
\end{equation}
Recall the definition
\begin{align*}
    \psi(\mathbf{k}_r;s)&=\frac{1}{\Gamma(s)}\int_0^\infty t^{s-1} \frac{{\rm A}(\mathbf{k}_r;\tanh t/2)}{\sinh(t)} dt,
\end{align*}
and we substitute equation (\ref{eq-2-3}) into the above equation and apply Proposition \ref{prop-2-1} to obtain the desired formula for $\psi(\mathbf{k}_r;s)$. \hfill$\square$

\begin{re} A result similar to Theorem \ref{thm-psi} for Arakawa-Kaneko zeta values can be found in \cite[Thm. 3.6]{KT2018}.

\end{re}

\subsection{Some explicit forms of Arakawa-Kaneko and Kaneko-Tsumura zeta functions}

Theorem \ref{thm-Ath} and Theorem \ref{thm-psi} can be written explicitly for some special arguments. In this section, we obtain some explicit forms of Theorem \ref{thm-Ath} and Theorem \ref{thm-psi}.

Let us consider the integral representation of the multiple polylogarithm ${\mathrm{Li}}_{{\bf k}_r}(z)$ and the level two multiple polylogarithm ${\rm A}({\bf k}_r;z)$ as follows.
\begin{align*}
&{\mathrm{Li}}_{{\bf k}_r}\left( z \right)=\int_{0<t_1<t_2\cdots<t_{k}<z}\frac{dt_1}{1-t_1}\underbrace{\frac{dt_2}{t_2}\cdots\frac{dt_{k_1}}{t_{k_1}}}_{(k_1-1)  \rm{-times}}\cdots \frac{dt_{k-k_r+1}}{1-t_{k-k_r+1}}\underbrace{\frac{dt_{k-k_r+2}}{t_{k-k_r+2}}\cdots\frac{dt_{k}}{t_{k}}}_{(k_r-1) \rm{-times}}
\end{align*}
and
\begin{align*}
    {\rm A}({\bf k}_r;z)
    =\int_{0<t_1<t_2\cdots<t_{k}<z}\frac{2dt_1}{1-t_1^2}&\underbrace{\frac{dt_2}{t_2}\cdots\frac{dt_{k_1}}{t_{k_1}}}_{(k_1-1)  \rm{-times}}\cdots \frac{2dt_{k-k_r+1}}{1-t_{k-k_r+1}^2}\underbrace{\frac{dt_{k-k_r+2}}{t_{k-k_r+2}}\cdots\frac{dt_{k}}{t_{k}}}_{(k_r-1) \rm{-times}},
\end{align*}
where ${\bf k}_r=({k_1}, \cdots ,k_{r-1},{k_r})$.
By using the formula
\begin{align}\label{b5}
  \int_{a<t_1\cdots<t_r<b}\underbrace{\frac{dt_1}{t_1}\cdots\frac{dt_r}{t_r}}_{r\text{-times}} =\frac{1}{r!}\left(\log \frac{b}{a}\right)^r,
\end{align}
we can write above integral expressions as
\begin{align}\label{2.6}
    {\mathrm{Li}}_{{\bf k}_r}(z)
    =\frac{1}{\prod_{i=1}^r(k_i-1)!}\int_{0<t_1<t_2\cdots<t_r<z}\frac{dt_1}{1-t_1}&\left(\log\frac{t_2}{t_1}\right)^{k_1-1}\cdots \frac{dt_r}{1-t_r}\left(\log\frac{z}{t_r}\right)^{k_r-1}
\end{align}
and
\begin{align}\label{2.7}
    {\rm A}({\bf k}_r;z)
    =\frac{2^r}{\prod_{i=1}^r(k_i-1)!}\int_{0<t_1<t_2\cdots<t_r<z}\frac{dt_1}{1-t_1^2}&\left(\log\frac{t_2}{t_1}\right)^{k_1-1}\cdots \frac{dt_r}{1-t_r^2}\left(\log\frac{z}{t_r}\right)^{k_r-1}
\end{align}
respectively.

Then, applying the changes of variables $t_j \mapsto 1-t_{r+1-j}$ and $t_j \mapsto \frac{1-t_{r+1-j}}{1+t_{r+1-j}}\quad (j=1,2,\ldots,r)$ to (\ref{2.6}) and (\ref{2.7}), respectively, we obtain
\begin{align}\label{b8}
{\mathrm{Li}}_{{\bf k}_r}\left( z \right)=\left\{\prod\limits_{j=1}^r\frac{(-1)^{k_j-1}}{(k_j-1)!}\right\}\int\nolimits_{E_r(z)} &\left\{\prod\limits_{j=1}^{r-1}\frac{\log^{k_j-1}\left(\frac{1-t_{r+1-j}}{1-t_{r-j}}\right)}{t_{r+1-j}}dt_{r+1-j}\right\}\nonumber\\ &\times\frac{\log^{k_r-1}\left(\frac{1-t_1}{z}\right)}{t_1}dt_1,
\end{align}
and
\begin{align}\label{b9}
&{\mathrm{A}}({\bf k}_r; z )=\left\{\prod\limits_{j=1}^r\frac{(-1)^{k_j-1}}{(k_j-1)!}\right\}\nonumber\\&\times\int\nolimits_{F_r(z)} \left\{\prod\limits_{j=1}^{r-1}\frac{\log^{k_j-1}\left(\frac{(1-t_{r+1-j})(1+t_{r-j})}{(1+t_{r+1-j})(1-t_{r-j})}\right)}{t_{r+1-j}}dt_{r+1-j}\right\}\frac{\log^{k_r-1}\left(\frac{1-t_1}{(1+t_1)z}\right)}{t_1}dt_1,
\end{align}
where
$$E_{r}(z):=\{(t_1,\ldots,t_r)\mid 1-z<t_1<\cdots<t_r<1\},$$
$$F_{r}(z):=\left\{(t_1,\ldots,t_r)\mid \frac{1-z}{1+z}<t_1<\cdots<t_r<1\right\}.$$
For the convenience (later use), let $E'_{r}(z):=E_{r}(1-z)$ and $F'_{r}(z):=F_{r}\left(\frac{1-z}{1+z}\right)$.

Let us consider the following lemma which will be needed in proving our main results under this section.
\begin{lem}\label{e3}
For integers $m\ge 0$ and $n>0$, we get
\begin{align}\label{e4}
\int\limits_{0}^z \log^m(t)\log^n\left(\frac{1-t}{1+t}\right)\frac{dt}{t}=(-1)^nn!\sum\limits_{l=0}^m l!\binom{m}{l} (-1)^l (\log(z))^{m-l} {\rm A}(\{1\}_{n-1},l+2;z).
\end{align}
In particular,
\begin{align}\label{e5}
\int\limits_{0}^1 \log^m(t)\log^n\left(\frac{1-t}{1+t}\right)\frac{dt}{t}&=(-1)^{n+m}n!m!T(\{1\}_{n-1},m+2).
\end{align}
\end{lem}
\pf
From Lemma (\ref{e2}), we have
\begin{align*}
\int\limits_{0}^z \log^m(t)\log^n\left(\frac{1-t}{1+t}\right)\frac{dt}{t}&=(-1)^n n!\int\limits_{0}^z \frac{\log^m(t){\rm A}(\{1\}_n;t)}{t}dt\\
&=(-1)^n n!\sum\limits_{l=0}^m l!\binom{m}{l} (-1)^l (\log(z))^{m-l} {\rm A}(\{1\}_{n-1},l+2;z).
\end{align*}

By setting $z=1$ in the above equation, we get
\begin{align*}
\int\limits_{0}^1 \log^m(t)\log^n\left(\frac{1-t}{1+t}\right)\frac{dt}{t}&=(-1)^{n+m}n!m!{\rm A}(\{1\}_{n-1},m+2;1)\nonumber\\
&=(-1)^{n+m}n!m!T(\{1\}_{n-1},m+2).
\end{align*}
This completes the proof of the lemma. \hfill$\square$

\begin{thm}\label{thm2.4} For any positive integers $j$ and $r$ with $j\leq r$,
\begin{align*}
{\rm A}\left(\{1\}_{j-1},2,\{1\}_{r-j};\frac{1-z}{1+z}\right)=&\sum\limits_{i=0}^{r-j}(-1)^{i}\binom{i+j}{i} T(i+j+1){\rm A}\left(\{1\}_{r-j-i};\frac{1-z}{1+z}\right)\nonumber\\
&+(-1)^{r-j-1}\sum\limits_{l=r-j}^{r}\binom{l}{r-j} {\rm A}\left(\{1\}_{r-l};\frac{1-z}{1+z}\right) {\rm A}(l+1;z).
\end{align*}
\end{thm}
\pf The proof of Theorem \ref{thm2.4} is similar to the proof of Theorem \ref{thm2.3}. Letting $k_1=\cdots=k_{j-1}=1,k_j=2,k_{j+1}=\cdots=k_r=1$ and replacing $z$ by $\frac{1-z}{1+z}$ in (\ref{b9}), using (\ref{b5}), we have
\begin{align}\label{b16}
&\A\left(\{1\}_{j-1},2,\{1\}_{r-j};\frac{1-z}{1+z}\right)\nonumber\\=&-\int\nolimits_{F'_r(z)}\frac{\log\left(\frac{(1-t_{r+1-j})(1+t_{r-j})}{(1+t_{r+1-j})(1-t_{r-j})}\right)}{t_1\cdots t_r}dt_1\cdots dt_r\nonumber\\
=&\sum\limits_{i=0}^{r-j} \frac{(-1)^{r-i}}{(r-j-i)!} \frac{i+j}{i!j!} \log^{r-j-i}(z) \int\limits_{z}^1 \frac{\log^{i+j-1}(t)\log\left(\frac{1-t}{1+t}\right)}{t}dt.
\end{align}
Substituting (\ref{e4}) and (\ref{e5}) with $n=1$ into (\ref{b16}). Then, we get
\begin{align}
\A&\left(\{1\}_{j-1},2,\{1\}_{r-j};\frac{1-z}{1+z}\right)\nonumber\\
&=\sum\limits_{i=0}^{r-j} \frac{(-1)^{r+j}}{(r-j-i)!} \binom{i+j}{i} \log^{r-j-i}(z) T(i+j+1)\nonumber\\
&\;\;\;\;\;\;\;\;\;\;\;\;\;\;\;+\sum\limits_{i=0}^{r-j}\sum\limits_{l=0}^{i+j-1} \frac{(-1)^{r-i+l}}{(r-l-1)!} \binom{i+j}{i}\binom{r-l-1}{r-j-i} \log^{r-l-1}(z) {\rm A}(l+2;z).
\end{align}
By substituting Lemma \ref{e2} into the above equation, we get

\begin{align}\label{eq:4.12}
\A&\left(\{1\}_{j-1},2,\{1\}_{r-j};\frac{1-z}{1+z}\right)\nonumber\\
&=\sum\limits_{i=0}^{r-j} (-1)^{i} \binom{i+j}{i}T(i+j+1) \A\left(\{1\}_{r-j-i};\frac{1-z}{1+z}\right) \nonumber
\\
&\quad+\sum\limits_{i=0}^{r-j}\sum\limits_{l=0}^{i+j-1} (-1)^{i+1} \binom{i+j}{i}\binom{r-l-1}{r-j-i} \A\left(\{1\}_{r-l-1};\frac{1-z}{1+z}\right) \A(l+2;z).
\end{align}
In order to obtain the desired formula, let us simplify the last term as follows.
\begin{align*}
   \sum\limits_{i=0}^{r-j}\sum\limits_{l=0}^{i+j-1} &(-1)^{i+1} \binom{i+j}{i}\binom{r-l-1}{r-j-i} \A\left(\{1\}_{r-l-1};\frac{1-z}{1+z}\right) \A(l+2;z)
   \\
&=\sum\limits_{i=0}^{r-j}\sum\limits_{l=1}^{i+j} (-1)^{i+1} \binom{i+j}{i}\binom{r-l}{r-j-i} \A\left(\{1\}_{r-l};\frac{1-z}{1+z}\right) \A(l+1;z)
\\
&=\sum\limits_{l=1}^{r}\sum\limits_{i=0}^{r-j} (-1)^{i+1} \binom{i+j}{i}\binom{r-l}{r-j-i} \A\left(\{1\}_{r-l};\frac{1-z}{1+z}\right) \A(l+1;z).
\end{align*}
Here,
\begin{align*}
\binom{r-l}{r-j-i}= \binom{r-l}{i+j-l}.
\end{align*}
By using the binomial identity 176 in \cite{michael}, we get
\begin{align*}
\sum\limits_{i=0}^{r-j} (-1)^{i+1} \binom{i+j}{i}\binom{r-l}{i+j-l}= \left\{ {\begin{array}{*{20}{c}} 0
   \;\;\;\;\;\;\;\;\;\;\;\;\;\;\;{\ \ (j < r-l),}  \\
   {(-1)^{r-j-1}\binom{l}{r-j}\;\;\;(r-l\leq j\leq r).}  \\
\end{array} } \right.
\end{align*}
Substituting this into the above equation, we get
\begin{align*}
   \sum\limits_{i=0}^{r-j}\sum\limits_{l=0}^{i+j-1} &(-1)^{i+1} \binom{i+j}{i}\binom{r-l-1}{r-j-i} \A\left(\{1\}_{r-l-1};\frac{1-z}{1+z}\right) \A(l+2;z)
   \\
   &=\sum\limits_{l=r-j}^{r}(-1)^{r-j-1}\binom{l}{r-j} \A\left(\{1\}_{r-j};\frac{1-z}{1+z}\right) \A(l+1;z).
\end{align*}
By substituting this in equation (\ref{eq:4.12}), we get the desired results.\hfill$\square$

Similarly, using (\ref{b5}) and (\ref{b8}) we can obtain the following formula for multiple polylogarithm function.
\begin{thm}\label{thm2.3} For any positive integers $j$ and $r$ with $j\leq r$,
\begin{align}\label{b11}
{\rm Li}_{\{1\}_{j-1},2,\{1\}_{r-j}}(1-z)=&\sum\limits_{i=0}^{r-j}(-1)^{i}\binom{i+j}{i} \zeta(i+j+1){\rm Li}_{\{1\}_{r-j-i}}(1-z)\nonumber\\
&+(-1)^{r-j-1}\sum\limits_{l=r-j}^{r}\binom{l}{r-j} {\rm Li}_{\{1\}_{r-l}}(1-z) {\rm Li}_{l+1}(z).
\end{align}
\end{thm}
\pf The proof of Theorem \ref{thm2.3} is similar to the proof of Theorem \ref{thm2.4} and is thus omitted. We leave the detail to the
interested reader.\hfill$\square$

In order to prove the next result we consider the following shuffle product identity.
\begin{lem}\label{lm-4.14}
For the integers $m, n \ge1$, we have
\begin{align*}
\sum\limits_{j=1}^{m}(-1)^j (y^{m-j} \shuffle y^j x^n)=-\sum_{\substack{\alpha_1+\cdots+\alpha_m=m+n, \forall \alpha_i \ge 1 }} y x^{\alpha_1} \cdots yx^{\alpha_{m-1}}yx^{\alpha_{m}}.
\end{align*}
\end{lem}
\pf
Consider the left hand-side of the above equation.
\begin{align}\label{4.17}
\sum\limits_{j=1}^{m}(-1)^j& y^{m-j} \shuffle y^j x^n \nonumber\\
&=(-1)^m y^mx^n+\sum\limits_{j=1}^{m-1} (-1)^j \left( y(y^{m-j-1} \shuffle y^j x^n)+y(y^{m-j} \shuffle y^{j-1} x^n) \right) \nonumber\\
&=(-1)^m y^mx^n+\sum\limits_{j=2}^{m} (-1)^{j-1}  y(y^{m-j} \shuffle y^{j-1} x^n)+ \sum\limits_{j=1}^{m-1}(-1)^j y(y^{m-j} \shuffle y^{j-1} x^n)  \nonumber\\
&=(-1)^m y^mx^n+(-1)^{m-1} y^mx^n-y(y^{m-1} \shuffle x^n) \nonumber \\
&=-y(y^{m-1} \shuffle x^n) .
\end{align}
By using the shuffle product formula, we obtain the desired result.\hfill$\square$

By using the above mentioned lemma, we can obtain the following formula for ${\rm A}(\mathbf{k}_r;z)$.
\begin{pro}\label{2.8} For $n,m\in \N$, we have
\begin{align}\label{b21}
\sum\limits_{j=1}^{m}(-1)^j{\rm A}(\{1\}_{m-j};z){\rm A}(\{1\}_{j-1},n+1;z)=-\sum_{\substack{|{\bf{k}}'|=m+n, \forall k_i \ge 1 \\d({\bf{k}}')=m}}{\rm A}({\bf{k}}';z).
\end{align}
\end{pro}
\pf We know that ${\rm A}(\mathbf{k}_r;z)$ satisfies the shuffle relation. Using Lemma \ref{lm-4.14}, we can obtain the desired results. \hfill$\square$

\begin{thm}\label{thm2.5} For any positive integers $r$ and $k$,
\begin{align*}
{\rm A}\left(\{1\}_{r-1},k;\frac{1-z}{1+z}\right)
&=\sum\limits_{j=0}^{k-2} (-1)^{k-j}T(\{1\}_{j},r+1){\rm A}(\{1\}_{k-2-j};z)\nonumber\\
&\quad +(-1)^{k-1} \sum_{\substack{a_1+\cdots+a_k=r \\ \forall a_j \ge 0}}{\rm A}\left(\{1\}_{a_k};\frac{1-z}{1+z}\right){\rm A}(a_1+1,\cdots,a_{k-1}+1;z) .
\end{align*}
\end{thm}
\pf
Set $k_1=\cdots=k_{r-1}=1,k_r=k$ in (\ref{b9})  and replacing $z$ by $\frac{1-z}{1+z}$. Then, we get
\begin{align*}
{\rm A}\left(\{1\}_{r-1},k;\frac{1-z}{1+z}\right)&=\frac{1}{(k-1)!}\int\nolimits_{F'_r(z)}\log^{k-1}\frac{(1-z)(1+t_1)}{(1+z)(1-t_1)}\frac{dt_1}{t_1}\cdots \frac{dt_r}{t_r}\nonumber\\
&=\sum\limits_{j=1}^{k-1} \frac{(-1)^j}{(k-1-j)!j!}\log^{k-1-j}\left(\frac{1-z}{1+z}\right)\int\nolimits_{F'_r(z)} \log^j\left(\frac{1-t_1}{1+t_1}\right)\frac{dt_1}{t_1}\cdots \frac{dt_r}{t_r}\nonumber
\\
&\;\;\;\;\;\;\quad+\frac{1}{(k-1)!}\log^{k-1}\left(\frac{1-z}{1+z}\right)\int\limits_{F'_r(z)} \frac{dt_1}{t_1}\cdots \frac{dt_r}{t_r}.
\end{align*}
By using (\ref{b5}), we get
\begin{align}\label{b20}
{\rm A}\left(\{1\}_{r-1},k;\frac{1-z}{1+z}\right)&=\sum\limits_{j=1}^{k-1} \frac{(-1)^{r-1+j}}{(k-1-j)!j!(r-1)!}\log^{k-1-j}\left(\frac{1-z}{1+z}\right)\int\nolimits_{z}^1 \log^{r-1}(t)\log^j\left(\frac{1-t}{1+t}\right)\frac{dt}{t}\nonumber
\\
&\;\;\;\;\;\;\quad+\frac{(-1)^r}{(k-1)!r!}\log^{k-1}\left(\frac{1-z}{1+z}\right)\log^r(z).
\end{align}
Substitute (\ref{e4}) and (\ref{e5}) into above equation. Then, we get
\begin{align}
{\rm A}&\left(\{1\}_{r-1},k;\frac{1-z}{1+z}\right)\nonumber\\
&=\sum\limits_{j=1}^{k-1} \frac{1}{(k-1-j)!} \log^{k-1-j}\left(\frac{1-z}{1+z}\right) T(\{1\}_{j-1},r+1)\nonumber\\
&\;\;\;\;\;\;\;\;\;\;\;+\sum\limits_{j=1}^{k-1}\sum\limits_{l=0}^{r-1} \frac{(-1)^{r-l}}{(k-1-j)!(r-1-l)!} \log^{k-1-j}\left(\frac{1-z}{1+z}\right)\log^{r-l-1}(z) {\rm A}(\{1\}_{j-1},l+2;z)\nonumber\\
&\;\;\;\;\;\;\;\;\;\;\;\quad+\frac{(-1)^r}{(k-1)!r!}\log^{k-1}\left(\frac{1-z}{1+z}\right)\log^r(z).
\end{align}
By substituting Lemma \ref{e2} into the above equation, we get
\begin{align*}
{\rm A}&\left(\{1\}_{r-1},k;\frac{1-z}{1+z}\right)\\
&=(-1)^{k-1}{\rm A}(\{1\}_{k-1};z){\rm A}\left(\{1\}_r;\frac{1-z}{1+z}\right)+\sum\limits_{j=1}^{k-1} (-1)^{k-1-j} {\rm A}(\{1\}_{k-1-j};z) T(\{1\}_{j-1},r+1)
\end{align*}

\begin{align*}
\;\;\;\;\;\;\;\;\;\;\;+\sum\limits_{j=1}^{k-1}\sum\limits_{i=0}^{r-1} (-1)^{k-j}{\rm A}(\{1\}_{k-1-j};z) {\rm A}\left(\{1\}_{r-i-1};\frac{1-z}{1+z}\right) {\rm A}(\{1\}_{j-1},i+2;z).
\end{align*}

This can be written as
\begin{align}\label{eq:4.16}
{\rm A}&\left(\{1\}_{r-1},k;\frac{1-z}{1+z}\right)\nonumber\\
&=\sum\limits_{j=1}^{k-1} (-1)^{k-1-j}T(\{1\}_{j-1},r+1){\rm A}(\{1\}_{k-1-j};z)\nonumber\\
&\;\;\;\;\;\;\;\;+(-1)^k\sum\limits_{i=0}^{r-1}{\rm A}\left(\{1\}_i;\frac{1-z}{1+z}\right)\Bigg(\sum\limits_{j=1}^{k-1}(-1)^j{\rm A}(\{1\}_{k-1-j};z){\rm A}(\{1\}_{j-1},r+1-i;z) \Bigg)\nonumber\\
&\;\;\;\;\;\;\;\;\;\;\;\;+(-1)^{k-1} {\rm A}\left(\{1\}_r;\frac{1-z}{1+z}\right){\rm A}(\{1\}_{k-1};z).
\end{align}
Setting $m=k-1$ and $n=r-j$ in Proposition \ref{2.8}, we can write the inner summation of the second term of equation (\ref{eq:4.16}) as below.
\begin{align*}
\sum\limits_{j=1}^{k-1}(-1)^j{\rm A}(\{1\}_{k-1-j};z){\rm A}(\{1\}_{j-1},r+1-i;z)=-\sum_{\substack{|{\bf{k}}'|=k-1+r-i, \forall k_i \ge 1 \\d({\bf{k}}')=k-1}}{\rm A}({\bf{k}}';z).
\end{align*}
By substituting this into equation (\ref{eq:4.16}), we get
\begin{align}
{\rm A}\left(\{1\}_{r-1},k;\frac{1-z}{1+z}\right)
&=\sum\limits_{j=0}^{k-2} (-1)^{k-j}T(\{1\}_{j},r+1){\rm A}(\{1\}_{k-2-j};z)\nonumber\\
&\quad+(-1)^{k-1}\sum\limits_{i=0}^{r-1}{\rm A}\left(\{1\}_i;\frac{1-z}{1+z}\right)\sum_{\substack{|{\bf{k}}'|=k-1+r-i, \forall k_i \ge 1 \\d({\bf{k}}')=k-1}}{\rm A}({\bf{k}}';z)\nonumber\\
&\quad+(-1)^{k-1} {\rm A}\left(\{1\}_r;\frac{1-z}{1+z}\right){\rm A}(\{1\}_{k-1};z)\nonumber\\
&=\sum\limits_{j=0}^{k-2} (-1)^{k-j}T(\{1\}_{j},r+1){\rm A}(\{1\}_{k-2-j};z)\nonumber\\
&\quad+(-1)^{k-1}\sum\limits_{i=0}^{r}{\rm A}\left(\{1\}_i;\frac{1-z}{1+z}\right)\sum_{\substack{|{\bf{k}}'|=k-1+r-i, \forall k_i \ge 1 \\d({\bf{k}}')=k-1}}{\rm A}({\bf{k}}';z).
\end{align}
We can write the second term of the above equation as
\begin{align}
(-1)^{k-1} \sum_{\substack{a_1+\cdots+a_k=r \\ \forall a_j \ge 0}}{\rm A}\left(\{1\}_{a_k};\frac{1-z}{1+z}\right){\rm A}(a_1+1,\cdots,a_{k-1}+1;z).
\end{align}
From this we can obtain the desired result.\hfill$\square$

By setting $j=1$ in Theorems \ref{thm2.4} and \ref{thm2.3}, we obtain the following corollaries respectively.
\begin{cor}(\cite{KTB2018}) For positive integer $r$,
\begin{align}\label{b25}
\A\left(2,\{1\}_{r-1};\frac{1-z}{1+z}\right)&=(-1)^rr\A(r+1;z)-(-1)^r\A(r;z)\log(z)\nonumber\\
&\quad-(-1)^r\sum\limits_{i=0}^{r-1}\frac{i+1}{(r-1-i)!}T(i+2)\log^{r-1-i}(z).
\end{align}
\end{cor}

\begin{cor}(\cite{KTA2018}) For positive integer $r$,
\begin{align}\label{b24}
{\rm Li}_{2,\{1\}_{r-1}}(1-z)&=(-1)^rr{\rm Li}_{r+1}(z)-(-1)^r{\rm Li}_r(z)\log(z)\nonumber\\
&\quad-(-1)^r\sum\limits_{i=0}^{r-1}\frac{i+1}{(r-1-i)!}\zeta(i+2)\log^{r-1-i}(z).
\end{align}
\end{cor}

Now, accordingly, we obtain explicit formulas for $\psi({\bf k}_r;s)$ and $\xi({\bf k}_r;s)$.

\begin{thm}\label{th:4.11} For positive integers $j,r$ and $\Re(s)>1$ with $j\leq r$,
\begin{align*}
\psi&(\{1\}_{j-1},2,\{1\}_{r-j};s)\\
&=\sum\limits_{i=0}^{r-j} (-1)^i \binom{i+j}{i} \binom{s+r-i-j-1}{r-i-j}T(i+j+1)T(s+r-i-j)\nonumber\\
&\;\;\;\;\;\;\;\;\;\;\;\;+(-1)^{r-j-1}\sum\limits_{l=r-j}^{r} \binom{l}{r-j}\binom{s+r-l-1}{r-l}T(l+1,s+r-l).
\end{align*}
\end{thm}
\pf
Let us consider Theorem \ref{thm2.5}.
\begin{align}\label{4.15}
{\rm A}\left(\{1\}_{j-1},2,\{1\}_{r-j};\frac{1-z}{1+z}\right)=&\sum\limits_{i=0}^{r-j}(-1)^{i}\binom{i+j}{i} T(i+j+1)\A\left(\{1\}_{r-j-i};\frac{1-z}{1+z}\right)\nonumber\\
&+(-1)^{r-j-1}\sum\limits_{l=r-j}^{r}\binom{l}{r-j} \A\left(\{1\}_{r-l};\frac{1-z}{1+z}\right) \A(l+1;z).
\end{align}
Now we can see that the right side the above equation is in the form of Theorem \ref{thm-Ath}. We can write the each term of the equation (\ref{4.15}) in the form of Theorem \ref{thm-psi}. This readily gives the desired result. \hfill$\square$

Similarly, we can obtain the formula for $\xi({\bf k}_r;s)$ as follows.
\begin{thm}\label{thm2.12} For positive integers $j,r$ and $\Re(s)>1$ with $j\leq r$,
\begin{align}\label{b29}
\xi(\{1\}_{j-1},2,\{1\}_{r-j};s)=&\sum\limits_{i=0}^{r-j} (-1)^i \binom{i+j}{i} \binom{s+r-i-j-1}{r-i-j}\zeta(i+j+1)\zeta(s+r-i-j)\nonumber\\
&+(-1)^{r-j-1}\sum\limits_{l=r-j}^{r} \binom{l}{r-j}\binom{s+r-l-2}{r-l-1}\zeta(l+2,s+r-l-1).
\end{align}
\end{thm}

Note that the explicit formulas of $\xi(2,\{1\}_{r-1};s)$ and $\psi(2,\{1\}_{r-1};s)$ were given by Kaneko and Tsumura in \cite{KT2018} and \cite{KTA2018}, respectively.

\section{Duality relation for Kaneko-Tsumura $\psi$-values}
In \cite{KTA2018}, Kaneko and Tsumura gave the following duality formula
\begin{align}\label{c1}
&\psi(\{1\}_{r-1},k;m+1)+(-1)^k\psi(\{1\}_{m-1},k;r+1)\nonumber\\
&=\sum\limits_{j=0}^{k-2} (-1)^j T(\{1\}_{r-1},k-j) T(\{1\}_j,m+1).
\end{align}
In this section, we will give a general duality formula for Kaneko-Tsumura $\psi$-values.

\subsection{Multiple $T$-harmonic sums and multiple $S$-harmonic sums}

Firstly, we define the multiple $T$-harmonic sums and the multiple $S$-harmonic sums, which can be regarded as two variants of the classical multiple harmonic sums of level two.
For indexes ${\bfk_r}:= (k_1,\ldots, k_r)\in \mathbb{Z}_{\ge 1}^r$ and ${\bfl_s}:=(l_1,l_2,\ldots,l_s)\in \mathbb{Z}_{\ge 1}^r$, and any positive integers $m$ and $p$, let
\begin{align*}
\bfk_{2m-1}:=(k_1,k_2,\ldots,k_{2m-1}),\quad \bfk_{2m}:=(k_1,k_2,\ldots,k_{2m})
\end{align*}
and
\begin{align*}
\bfl_{2p-1}:=(l_1,l_2,\ldots,l_{2p-1}),\quad \bfl_{2p}:=(l_1,l_2,\ldots,l_{2p}).
\end{align*}

For positive integers $n_1,n_2,\ldots,n_{r}$ and $n$, if $r=2m-1$ is an odd, we define
\begin{align*}
&D_n({{\bf n}_{2m-1}}):=\left\{(n_1,n_2,\ldots,n_{2m-1},n)\mid 0<n_1\leq n_2 <\cdots\leq n_{2m-2}<n_{2m-1}\leq n \right\},\ (n\geq m)\\
&E_n({{\bf n}_{2m-1}}):=\left\{(n_1,n_2,\ldots,n_{2m-1},n)\mid 1\leq n_1<n_2\leq \cdots< n_{2m-2}\leq n_{2m-1}< n \right\},\ (n>m)
\end{align*}
and when $r=2m$ is even, we define
\begin{align*}
&D_n({{\bf n}_{2m}}):=\left\{(n_1,n_2,\ldots,n_{2m},n)\mid 0<n_1\leq n_2 <\cdots\leq n_{2m-2}<n_{2m-1}\leq n_{2m}<n \right\},\ (n> m)\\
&E_n({{\bf n}_{2m}}):=\left\{(n_1,n_2,\ldots,n_{2m},n)\mid 1\leq n_1<n_2\leq \cdots< n_{2m-2}\leq n_{2m-1}< n_{2m}\leq n \right\},\ (n>m).
\end{align*}

\begin{defn}\emph{(\cite[Def. 1.1]{XZ2020})}\label{def1} For positive integer $m$, the multiple $T$-harmonic sums ({\rm MTHSs} for short) and multiple $S$-harmonic sums ({\rm MSHSs} for short) are defined by
\begin{align}
&T_n({\bfk_{2m-1}}):= \sum_{D_n({{\bf n}_{2m-1}})} \frac{2^{2m-1}}{(\prod_{j=1}^{m-1} (2n_{2j-1}-1)^{k_{2j-1}}(2n_{2j})^{k_{2j}})(2n_{2m-1}-1)^{k_{2m-1}}},\label{MOT}\\
&T_n({\bfk_{2m}}):= \sum_{D_n({{\bf n}_{2m}})} \frac{2^{2m}}{\prod_{j=1}^{m} (2n_{2j-1}-1)^{k_{2j-1}}(2n_{2j})^{k_{2j}}},\label{MET}\\
&S_n({\bfk_{2m-1}}):= \sum_{E_n({{\bf n}_{2m-1}})} \frac{2^{2m-1}}{(\prod_{j=1}^{m-1} (2n_{2j-1})^{k_{2j-1}}(2n_{2j}-1)^{k_{2j}})(2n_{2m-1})^{k_{2m-1}}},\label{MOS}\\
&S_n({\bfk_{2m}}):= \sum_{E_n({{\bf n}_{2m}})} \frac{2^{2m}}{\prod_{j=1}^{m} (2n_{2j-1})^{k_{2j-1}}(2n_{2j}-1)^{k_{2j}}},\label{MES}
\end{align}
where $T_n({\bfk_{2m-1}}):=0$ if $n<m$, and $T_n({\bfk_{2m}})=S_n({\bfk_{2m-1}})=S_n({\bfk_{2m}}):=0$ if $n\leq m$. Moreover, for convenience we let $T_n(\emptyset)=S_n(\emptyset):=1$. We call \eqref{MOT} and \eqref{MET} are multiple $T$-harmonic sums, and call \eqref{MOS} and \eqref{MES} are multiple $S$-harmonic sums.
\end{defn}

Clearly, according to the definitions of MTHSs and MSHSs, we have the following relations
\begin{alignat*}{3}
&T_n({\bfk_{2m}})=2\sum_{j=1}^{n-1} \frac{T_j({\bfk_{2m-1}})}{(2j)^{k_{2m}}}, \qquad &
&T_n({\bfk_{2m-1}})=2\sum_{j=1}^{n} \frac{T_j({\bfk_{2m-2}})}{(2j-1)^{k_{2m-1}}},\\
&S_n({\bfk_{2m}})=2\sum_{j=1}^{n} \frac{S_j({\bfk_{2m-1}})}{(2j-1)^{k_{2m}}}, \qquad &
&S_n({\bfk_{2m-1}})=2\sum_{j=1}^{n-1} \frac{S_j({\bfk_{2m-2}})}{(2j)^{k_{2m-1}}}.
\end{alignat*}

In \cite{XZ2020}, the second author and Zhao used MTHSs and MSHSs to define the convoluted $T$-values $T({\bfk_{r}}\circledast {\bfl_{s}})$, which can be regarded as a $T$-variant of Kaneko-Yamamoto MZVs $\zeta({\bfk}_r\circledast {\bfl}_s)$ (see \cite{KY2018}).
\begin{defn} For positive integers $m$ and $p$, the \emph{convoluted $T$-values}
\begin{align}
&T({\bfk_{2m}}\circledast{\bfl_{2p}})=2\su \frac{T_n({\bfk_{2m-1}})T_n({\bfl_{2p-1}})}{(2n)^{k_{2m}+l_{2p}}},\\
&T({\bfk_{2m-1}}\circledast{\bfl_{2p-1}})=2\su \frac{T_n({\bfk_{2m-2}})T_n({\bfl_{2p-2}})}{(2n-1)^{k_{2m-1}+l_{2p-1}}},\\
&T({\bfk_{2m}}\circledast{\bfl_{2p-1}})=2\su \frac{T_n({\bfk_{2m-1}})S_n({\bfl_{2p-2}})}{(2n)^{k_{2m}+l_{2p-1}}},\\
&T({\bfk_{2m-1}}\circledast{\bfl_{2p}})=2\su \frac{T_n({\bfk_{2m-2}})S_n({\bfl_{2p-1}})}{(2n-1)^{k_{2m-1}+l_{2p}}}.
\end{align}
\end{defn}

Note that the MTVs are special cases of the convoluted $T$-values since
\begin{align*}
&T({\bfk_{r}}\circledast (1))=T((\bfk_r)_{+}),\quad T((1)\circledast {\bfl_{2p-1}})=T((\bfl_{2p-1})_{+}).
\end{align*}

In particular, in \cite{XZ2020}, the second author and Zhao shown the following theorem.
\begin{thm}(\cite[Thm. 4.5]{XZ2020}) For positive integers $l_1,l_2$ and index ${\bfk}_r=(k_1,k_2,\ldots,k_r)$, the convoluted T -values $T({\bfk}_r\circledast (l_1,l_2))$ can be expressed in terms of products of MTVs and Riemann zeta values.
\end{thm}

\subsection{Duality formula of Kaneko-Tsumura $\psi$-values}

Now, for index $\bfk_r=(k_1,k_2,\ldots,k_r)\in \N^r$ and $|\bfk_r|:=k_1+k_2+\cdots+k_r$, we adopt the following notations:
\begin{align*}
&\overrightarrow{\bfk_j}:=(k_1,k_2,\ldots ,k_j),\quad \overleftarrow{\bfk_j}:=(k_r,k_{r-1},\ldots,k_{r+1-j}),\\
&|\overrightarrow{\bfk_j}|:=k_1+k_2+\cdots+k_j,\quad |\overleftarrow{\bfk_j}|:=k_r+k_{r-1}+\cdots+k_{r+1-j}
\end{align*}
with $\overrightarrow{\bfk_0}=\overleftarrow{\bfk_0}:=\emptyset$ and $|\overrightarrow{\bfk_0}|=|\overleftarrow{\bfk_0}|:=0$.

\begin{thm}\label{thm3.1} For positive integers $p,q,r$ and index ${\bfk}=(k_1,\ldots,k_r)$ with $k_1,k_2,\ldots,k_r\in \N\setminus\{1\}$, let $\N^{+}_o$ be the set of positive odd numbers and $\N^{+}_e$ be the set of positive even numbers,
\begin{align}\label{c2}
&\psi\left(\{1\}_{q-1},\overrightarrow{\bfk}_r^{-};p+1\right)-(-1)^{|{\bfk}|}\psi\left(\{1\}_{p-1},\overleftarrow{\bfk}_r^{-};q+1\right)\nonumber\\
&=\sum\limits_{j=0}^{r-1} (-1)^{\mid\stackrel{\leftarrow}{{\bf k}}_j\mid}\sum\limits_{i=1}^{k_{r-j}-2} (-1)^{i-1} T\left(\{1\}_{p-1},\overleftarrow{\bfk}_j,i+1\right)T\left(\{1\}_{q-1},\overrightarrow{\bfk}_{r-j-1},k_{r-j}-i\right)\nonumber\\
&\quad+\sum\limits_{j=0}^{r-2} (-1)^{\mid\stackrel{\leftarrow}{{\bf k}}_{j+1}\mid}\left\{\begin{array}{l} T\left(\{1\}_{q-1},\overrightarrow{\bfk}_{r-j-1}\right)T\left(\Big(\{1\}_{p-1},\overleftarrow{\bfk}_{j+1}\Big)^{-}\circledast(1,1)\right)\\ - T\left(\{1\}_{p-1},\overleftarrow{\bfk}_{j+1}\right)T\left(\Big(\{1\}_{q-1},\overrightarrow{\bfk}_{r-j-1}\Big)^{-}\circledast(1,1)\right)\\
+2\delta_{p+j,q+r-j-2}\log(2)T\left(\{1\}_{p-1},\overleftarrow{\bfk}_{j+1}\right)T\left(\{1\}_{q-1},\overrightarrow{\bfk}_{r-j-1}\right)\end{array}\right\},
\end{align}
where ${\bfk_r^{-}}=(k_1,\ldots,k_{r-1},k_r-1)$ and
\begin{align*}
\delta_{r,s}= \left\{ {\begin{array}{*{20}{c}}\ \ 0 \quad(r,s\in \N^{+}_e)\quad\quad\quad\quad\\
   \ \ 0\quad(r,s\in \N^{+}_o)\quad\quad\quad\quad\\
   -1 \quad(r\in \N^{+}_e,\ s\in \N^{+}_o)\\
   \ \ \ 1 \quad(r\in \N^{+}_o,\ s\in \N^{+}_e).
\end{array} } \right.
\end{align*}
\end{thm}

It is clear that formula (\ref{c1}) is an immediate corollary of Theorem \ref{thm3.1} with $r=1$.

\subsection{Proof of duality formula}

In this section, we will prove the duality formula (\ref{c2}). We need the following a lemma.
\begin{lem}\label{lem3.5} Let sequences $A_n,B_n$ define the finite sums
 ${A_n} := \sum\limits_{k = 1}^n {{a_k}} ,\ {B_n} := \sum\limits_{k = 1}^n {{b_k}}\ ( {a_n},{b_n} =o(n^{-p}),\ {\mathop{\Re}\nolimits} \left( p \right) > 1 $ if $n\rightarrow \infty$)
 and $A = \mathop {\lim }\limits_{n \to \infty } {A_n},B = \mathop {\lim }\limits_{n \to \infty } {B_n}$, then
\begin{align*}
\su \left\{ \frac{A_nB}{n+\alpha}- \frac{AB_n}{n+\beta}\right\}=AB(\psi(\beta+1)-\psi(\alpha+1))+A\su b_nH_{n-1}(\beta)-B\su a_nH_{n-1}(\alpha),
\end{align*}
where $\alpha,\beta \notin \{-1,-2,-3,\ldots\},$ $\psi(\alpha+1)$ is digamma function, and $H_n(\alpha)$ is defined by
\[H_n(\alpha)=\sum_{k=1}^n\frac{1}{k+\alpha}.\]
It is clear that $T_n(1)=H_n(-1/2)$ and $S_n(1)=H_{n-1}(1)$.
\end{lem}
\pf The lemma is almost obvious. We leave the detail to the interested reader.\hfill$\square$

Now, we prove the following formula.
\begin{thm}\label{thm3.2} For positive integers $p,q,r$ and index $\bfk_r$ with $k_1,\ldots,k_r\in \N\setminus\{1\}$,
\begin{align}\label{c3}
&\psi\left(\{1\}_{q-1},\overrightarrow{\bfk}_r^{-};p+1\right)-(-1)^{|{\bfk}|}\psi\left(\{1\}_{p-1},\overleftarrow{\bfk}_r^{-};q+1\right)\nonumber\\
&=\sum\limits_{j=0}^{r-1} (-1)^{\mid\stackrel{\leftarrow}{{\bf k}}_j\mid}\sum\limits_{i=1}^{k_{r-j}-2} (-1)^{i-1} T\left(\{1\}_{p-1},\overleftarrow{\bfk}_j,i+1\right)T\left(\{1\}_{q-1},\overrightarrow{\bfk}_{r-j-1},k_{r-j}-i\right)\nonumber\\
&\quad+\sum\limits_{j=0}^{r-2} (-1)^{\mid\stackrel{\leftarrow}{{\bf k}}_{j+1}\mid} \lim_{x\rightarrow 1}\left\{\begin{array}{l} \A\left(\{1\}_{p-1},\overleftarrow{\bfk}_{j+1};x\right)\A\left(\{1\}_{q-1},\overrightarrow{\bfk}_{r-j-1},1;x\right)\\ -\A\left(\{1\}_{p-1},\overleftarrow{\bfk}_{j+1},1;x\right)\A\left(\{1\}_{q-1},\overrightarrow{\bfk}_{r-j-1};x\right)\end{array}\right\}.
\end{align}
\end{thm}
\pf We change variable $\tanh(t/2)=x$ in (\ref{a8}), and let $s=p+1$, then
\begin{align}
&\psi(k_1,k_2\ldots,k_r;p+1)=\frac{(-1)^p}{p!}\int\limits_{0}^1 \frac{\log^p\left(\frac{1-x}{1+x}\right)\A(k_1,k_2,\ldots,k_r;x)}{x}dx.\label{b35}
\end{align}
From (\ref{b4}) and (\ref{b35}), we can find that
\begin{align}\label{c3}
&\psi(\{1\}_{q-1},k_1,\ldots,k_{r-1},k_r-1;p+1)\nonumber
\\&=\int\limits_{0}^1 \frac{\A(\{1\}_p;u)\A(\{1\}_{q-1},k_1,\ldots,k_{r-1},k_r-1;u)}{u}du.
\end{align}
Then by using (\ref{b3}),
\begin{align*}
\frac{d}{dz}\A({{k_1}, \cdots ,k_{r-1},{k_r}}; z)= \left\{ {\begin{array}{*{20}{c}} \frac{1}{z}\A({{k_1}, \cdots ,{k_{r-1}},{k_r-1}};z)
   {\ \ (k_r\geq 2),}  \\
   {\frac{2}{1-z^2}\A({{k_1}, \cdots ,{k_{r-1}}};z)\;\;\;\ \ \ (k_r = 1).}  \\
\end{array} } \right.
\end{align*}
Hence, by integrating by parts, we give
\begin{align}\label{c4}
&\psi(\{1\}_{q-1},k_1,\ldots,k_{r-1},k_r-1;p+1)\nonumber
\\&= T(\{1\}_{p-1},2)T(\{1\}_{q-1},k_1,\ldots,k_{r-1},k_r-1)\nonumber\\&\quad-\int\limits_{0}^1 \frac{{\A}(\{1\}_{p-1},2;u){\A}(\{1\}_{q-1},k_1,\ldots,k_{r-1},k_r-2;u)}{u}du\nonumber\\
&=\cdots\nonumber\\
&=\sum\limits_{i=1}^{k_r-2}(-1)^{i-1}T(\{1\}_{p-1},i+1)T(\{1\}_{q-1},k_1,\ldots,k_{r-1},k_r-i)\nonumber\\
&\quad+(-1)^{k_r-2}\int\limits_{0}^1 \frac{{\A}(\{1\}_{p-1},k_r-1;u){\A}(\{1\}_{q-1},k_1,\ldots,k_{r-1},1;u)}{u}du\nonumber\\
&=\sum\limits_{i=1}^{k_r-2}(-1)^{i-1}T(\{1\}_{p-1},i+1)T(\{1\}_{q-1},k_1,\ldots,k_{r-1},k_r-i)\nonumber\\
&\quad+(-1)^{k_r}\lim_{x\rightarrow 1} \left\{\A(\{1\}_{p-1},k_r;x)\A(\{1\}_{q-1},k_1,\ldots,k_{r-1},1;x)\atop- 2\int\limits_{0}^1 \frac{\A(\{1\}_{p-1},k_r;u)\A(\{1\}_{q-1},k_1,\ldots,k_{r-1};u)}{1-u^2}du\right\}\nonumber\\
&=\sum\limits_{i=1}^{k_r-2}(-1)^{i-1}T(\{1\}_{p-1},i+1)T(\{1\}_{q-1},k_1,\ldots,k_{r-1},k_r-i)\nonumber\\
&\quad+(-1)^{k_r}\lim_{x\rightarrow 1} \left\{\A(\{1\}_{p-1},k_r;x)\A(\{1\}_{q-1},k_1,\ldots,k_{r-1},1;x)\atop-\A(\{1\}_{p-1},k_r,1;x)\A(\{1\}_{q-1},k_1,\ldots,k_{r-1};x)\right\}\nonumber\\
&\quad+(-1)^{k_r}\int\limits_{0}^1 \frac{\A(\{1\}_{p-1},k_r,1;u)\A(\{1\}_{q-1},k_1,\ldots,k_{r-2},k_{r-1}-1;u)}{u}du\nonumber\\
&=\cdots\nonumber\\
&=\sum\limits_{j=0}^{r-2} (-1)^{\mid\stackrel{\leftarrow}{{\bf k}}_j\mid}\sum\limits_{i=1}^{k_{r-j}-2} (-1)^{i-1} T(\{1\}_{p-1},k_r,k_{r-1},\ldots,k_{r+1-j},i+1)\nonumber\\&\quad\quad\quad\quad \quad\quad \quad\quad \quad\quad \times T(\{1\}_{q-1},k_1,k_{2},\ldots,k_{r-j-1},k_{r-j}-i)\nonumber\\
&\quad+\sum\limits_{j=0}^{r-2} (-1)^{\mid\stackrel{\leftarrow}{{\bf k}}_{j+1}\mid} \lim_{x\rightarrow 1}\left\{\begin{array}{l} \A(\{1\}_{p-1},k_r,k_{r-1},\ldots,k_{r-j};x)\A(\{1\}_{q-1},k_1,k_2,\ldots,k_{r-j-1},1;x)\\ -\A(\{1\}_{p-1},k_r,k_{r-1},\ldots,k_{r-j},1;x)\A(\{1\}_{q-1},k_1,k_2,\ldots,k_{r-j-1};x)\end{array}\right\}\nonumber\\
&\quad+(-1)^{\mid\stackrel{\leftarrow}{{\bf k}}_{r-1}\mid}\int\limits_{0}^1 \frac{\A(\{1\}_{p-1},k_r,k_{r-1},\ldots,k_2,1;u)\A(\{1\}_{q-1},k_1-1;u)}{u}du\nonumber\\
&=\sum\limits_{j=0}^{r-1} (-1)^{\mid\stackrel{\leftarrow}{{\bf k}}_j\mid}\sum\limits_{i=1}^{k_{r-j}-2} (-1)^{i-1} T(\{1\}_{p-1},k_r,k_{r-1},\ldots,k_{r+1-j},i+1)\nonumber\\&\quad\quad\quad\quad \quad\quad \quad\quad \quad\quad \times T(\{1\}_{q-1},k_1,k_{2},\ldots,k_{r-j-1},k_{r-j}-i)\nonumber\\
&\quad+\sum\limits_{j=0}^{r-2} (-1)^{\mid\stackrel{\leftarrow}{{\bf k}}_{j+1}\mid} \lim_{x\rightarrow 1}\left\{\begin{array}{l} \A(\{1\}_{p-1},k_r,k_{r-1},\ldots,k_{r-j};x)\A(\{1\}_{q-1},k_1,k_2,\ldots,k_{r-j-1},1;x)\\ -\A(\{1\}_{p-1},k_r,k_{r-1},\ldots,k_{r-j},1;x)\A(\{1\}_{q-1},k_1,k_2,\ldots,k_{r-j-1};x)\end{array}\right\}\nonumber\\
&\quad+(-1)^{\mid\stackrel{\leftarrow}{{\bf k}}_{r}\mid}\int\limits_{0}^1 \frac{\A(\{1\}_{p-1},k_r,k_{r-1},\ldots,k_2,k_1-1;u)\A(\{1\}_{q};u)}{u}du.
\end{align}
Hence, the formula (\ref{c3}) holds.\hfill$\square$

Next, we evaluate the limit in (\ref{c3}).

\begin{thm}\label{thm3.3} For indexes $\bfk_r$ and $\bfl_s$, and positive integers $r$ and $s$, we have
\begin{align}\label{c5}
&\lim_{x\rightarrow 1} \left\{\A(\bfk_r;x) \A(\bfl_s,1;x)-\A(\bfk_r,1;x) \A(\bfl_s;x)\right\}\nonumber\\
&=T(\bfl_s)T\left(\bfk_r^{-}\circledast(1,1)\right)-T(\bfk_r)T\left(\bfl_s^{-}\circledast(1,1)\right)+2\delta_{r,s}\log(2)T(\bfk_r)T(\bfl_s).
\end{align}
\end{thm}
\pf According to the definitions of ${\rm A}(k_1,k_2,\ldots,k_r;x)$ and MTHSs, by an elementary calculation, we can find that
\begin{align}\label{bc11}
&{\rm A}(\bfk_{2m-1};x)=2\su \frac{T_n(\bfk_{2m-2})}{(2n-1)^{k_{2m-1}}}x^{2n-1}
\end{align}
and
\begin{align}\label{bc12}
&{\rm A}(\bfk_{2m};x)=2\su \frac{T_n(\bfk_{2m-1})}{(2n)^{k_{2m}}}x^{2n}.
\end{align}

Hence, applying Lemma \ref{lem3.5},  by straightforward calculations, it is easy to see that if $r=2m-1$ and $s=2p-1$, then
\begin{align}\label{c6}
&\lim_{x\rightarrow 1} \left\{\A(\bfk_{2m-1};x) \A(\bfl_{2p-1},1;x)-\A(\bfk_{2m-1},1;x) \A(\bfl_{2p-1};x)\right\}\nonumber\\
&=\su \frac{T_n(\bfl_{2p-1})T(\bfk_{2m-1})-T_n(\bfk_{2m-1})T(\bfl_{2p-1})}{n}\nonumber\\
&=T(\bfl_{2p-1})2\su \frac{T_n(\bfk_{2m-2})S_n(1)}{(2n-1)^{k_{2m-1}}}-T(\bfk_{2m-1})2\su\frac{T_n(\bfl_{2p-2})S_n(1)}{(2n-1)^{l_{2p-1}}}\nonumber\\
&=T(\bfl_{2p-1})T\left(\bfk_{2m-1}^{-}\circledast(1,1) \right)-T(\bfk_{2m-1})T\left(\bfl_{2p-1}^{-}\circledast(1,1) \right),
\end{align}
if $r=2m-1$ and $s=2p$, then
\begin{align}\label{c7}
&\lim_{x\rightarrow 1} \left\{\A(\bfk_{2m-1};x) \A(\bfl_{2p},1;x)-\A(\bfk_{2m-1},1;x) \A(\bfl_{2p};x)\right\}\nonumber\\
&=\su \left\{\frac{T_n(\bfl_{2p})T(\bfk_{2m-1})}{n-1/2} -\frac{T_n(\bfk_{2m-1})T(\bfl_{2p})}{n} \right\}\nonumber\\
&=2\log(2)T(\bfl_{2p})T(\bfk_{2m-1})+T(\bfl_{2p})2\su \frac{T_n(\bfk_{2m-2})S_n(1)}{(2n-1)^{k_{2m-1}}}-T(\bfk_{2m-1})2\su\frac{T_n(\bfl_{2p-1})T_n(1)}{(2n)^{l_{2p}}}\nonumber\\
&=2\log(2)T(\bfl_{2p})+T(\bfl_{2p})T\left(\bfk_{2m-1}^{-}\circledast(1,1) \right)-T(\bfk_{2m-1})T\left(\bfl_{2p}^{-}\circledast(1,1) \right),
\end{align}
if $r=2m$ and $s=2p-1$, then
\begin{align}\label{c8}
&\lim_{x\rightarrow 1} \left\{\A(\bfk_{2m};x) \A(\bfl_{2p-1},1;x)-\A(\bfk_{2m},1;x) \A(\bfl_{2p-1};x)\right\}\nonumber\\
&=\su \left\{\frac{T_n(\bfl_{2p-1})T(\bfk_{2m})}{n}-\frac{T_n(\bfk_{2m})T(\bfl_{2p-1})}{n-1/2}  \right\}\nonumber\\
&=-2\log(2)T(\bfk_{2m})T(\bfl_{2p-1})+T(\bfl_{2p-1})2\su \frac{T_n(\bfk_{2m-1})T_n(1)}{(2n)^{k_{2m}}}-T(\bfk_{2m})2\su\frac{T_n(\bfl_{2p-2})S_n(1)}{(2n-1)^{l_{2p-1}}}\nonumber\\
&=-2\log(2)T(\bfk_{2m})T(\bfl_{2p-1})+T(\bfl_{2p-1})T\left(\bfk_{2m}^{-}\circledast(1,1) \right)-T(\bfk_{2m})T\left(\bfl_{2p-1}^{-}\circledast(1,1) \right),
\end{align}
if $r=2m$ and $s=2p$, then
\begin{align}\label{c9}
&\lim_{x\rightarrow 1} \left\{\A(\bfk_{2m};x) \A(\bfl_{2p},1;x)-\A(\bfk_{2m},1;x) \A(\bfl_{2p};x)\right\}\nonumber\\
&=\su \frac{T_n(\bfl_{2p})T(\bfk_{2m})-T_n(\bfk_{2m})T(\bfl_{2p})}{n-1/2}\nonumber\\
&=T(\bfl_{2p})2\su \frac{T_n(\bfk_{2m-2})T_n(1)}{(2n)^{k_{2m}}}-T(\bfk_{2m})2\su\frac{T_n(\bfl_{2p-1})T_n(1)}{(2n)^{l_{2p}}}\nonumber\\
&=T(\bfl_{2p})T\left(\bfk_{2m}^{-}\circledast(1,1) \right)-T(\bfk_{2m})T\left(\bfl_{2p}^{-}\circledast(1,1) \right).
\end{align}
Thus, combining (\ref{c6})-(\ref{c9}), we deduce the desired result.\hfill$\square$

Therefore, from Theorem \ref{thm3.2}, we can get the following formula
\begin{align}\label{c10}
&\lim_{x\rightarrow 1}\left\{\begin{array}{l} \A\left(\{1\}_{p-1},\overleftarrow{\bfk}_{j+1};x\right)\A\left(\{1\}_{q-1},\overrightarrow{\bfk}_{r-j-1},1;x\right)\\ -\A\left(\{1\}_{p-1},\overleftarrow{\bfk}_{j+1},1;x\right)\A\left(\{1\}_{q-1},\overrightarrow{\bfk}_{r-j-1};x\right)\end{array}\right\}\nonumber\\
&=T\left(\{1\}_{q-1},\overrightarrow{\bfk}_{r-j-1}\right)T\left(\Big(\{1\}_{p-1},\overleftarrow{\bfk}_{j+1}\Big)^{-}\circledast(1,1)\right)\nonumber\\&\quad - T\left(\{1\}_{p-1},\overleftarrow{\bfk}_{j+1}\right)T\left(\Big(\{1\}_{q-1},\overrightarrow{\bfk}_{r-j-1}\Big)^{-}\circledast(1,1)\right)\nonumber\\
&\quad+2\delta_{p+j,q+r-j-2}\log(2)T\left(\{1\}_{p-1},\overleftarrow{\bfk}_{j+1}\right)T\left(\{1\}_{q-1},\overrightarrow{\bfk}_{r-j-1}\right).
\end{align}
{\bf Proof of Theorem \ref{thm3.1}.} Substituting (\ref{c10}) into (\ref{c3}) yields the desired evaluation (\ref{c2}). \hfill$\square$

\section{Further discussion}\label{sec3}

{\bf Question}: Is there a function $y=y(z)$ that satisfies the following system for any $m\in\N$?
\begin{align*}
 \left\{ {\begin{array}{*{20}{c}} \frac{dy}{dz}+y^m=1,
    \\ \\
   \frac{d\log(y)}{dz}=\frac{m^2e^{-mz}}{1-e^{-m^2z}}, \\
\end{array} } \right.
\end{align*}
with the far-field boundary condition and the initial condition:
\begin{align*}
 \left\{ {\begin{array}{*{20}{c}} y(0)=0,
    \\ \\
   y(z)\rightarrow1\quad {\rm as} \quad z\rightarrow +\infty.\\
\end{array} } \right.
\end{align*}
Does the solution of this system exists?
In particular, if $m=1$, then $y(z)=1-e^{-z}$. If $m=2$, then $y(z)=\tanh(z)$. But for $m\geq 3$, does the solution of this system exists?

{\bf Next, we assume this solution $y=f_m(z)$ exists. Hence, $f_1(z)=1-e^{-z}$ and $f_2(z)=\tanh(z)$.}

\begin{defn} For $k_1,\ldots,k_r\in\N$, the multiple polylogarithm of level $m$ is defined by
\begin{align}\label{d1}
{\rm Ath}^{(m)}(k_1,k_2,\ldots,k_r;z): = \sum\limits_{1 \le {n_1} <  \cdots  < {n_r}\atop n_i\equiv i\ {\rm mod}\ m} {\frac{{{z^{{n_r}}}}}{{n_1^{{k_1}}n_2^{{k_2}} \cdots n_r^{{k_r}}}}},\quad z \in \left[ { - 1,1} \right).
\end{align}
\end{defn}
Note that if $z\in [0,+\infty)$, then $|f_m(z)|< 1$, we have
\begin{align}\label{d2}
{\rm Ath}^{(m)}(1;f_m(z))=\sum\limits_{n=1}^\infty \frac{f_m^{mn+1}(z)}{mn+1}=\sum\limits_{n=1}^\infty \int\limits_{0}^z f_m^{mn}(x)f'_m(x) dx =\int\limits_{0}^z\frac{f'_m(x)}{1-f_m^m(x)}dx=z.
\end{align}
Similar to (\ref{b3}), we can easily obtain the following.
\begin{lem} {\rm (i)} For $r,k_1,\ldots,k_r\in \N$,
\begin{align}\label{d3}
\frac{d}{dz}{\mathrm{Ath}^{(m)}}({{k_1}, \cdots ,k_{r-1},{k_r}}; z)= \left\{ {\begin{array}{*{20}{c}} \frac{1}{z} {\mathrm{Ath}^{(m)}}({{k_1}, \cdots ,{k_{r-1}},{k_r-1}};z)
   {\ \ (k_r\geq 2),}  \\
   {\frac{1}{1-z^m}{\mathrm{Ath}^{(m)}}({{k_1}, \cdots ,{k_{r-1}}};z)\;\;\;\ \ \ (k_r = 1).}  \\
\end{array} } \right.
\end{align}
{\rm (ii)} For $r\in\N$,
  \begin{align}\label{d4}
{\rm Ath}^{(m)}({\{1\}_r};z)=\frac{1}{r!}({\rm Ath}^{(m)}(1;z))^r.
\end{align}
\end{lem}
By (\ref{d3}), we obtain
\begin{align}\label{d5}
&{\mathrm{Ath}^{(m)}}({{k_1}, \cdots,k_{r-1} ,{k_r}};z)=\int\limits_{0}^z \underbrace{\frac{dt}{t}\cdots\frac{dt}{t}}_{k_r-1}\frac{dt}{1-t^m}\underbrace{\frac{dt}{t}\cdots\frac{dt}{t}}_{k_{r-1}-1}\frac{dt}{1-t^m}\cdots
\underbrace{\frac{dt}{t}\cdots\frac{dt}{t}}_{k_1-1}\frac{dt}{1-t^m}\nonumber\\
&=\left\{\prod\limits_{j=1}^r\frac{(-1)^{k_j-1}}{(k_j-1)!}\right\}\int\nolimits_{D_r(z)} \frac{\log^{k_1-1}\left(\frac{t_1}{t_2}\right)\cdots \log^{k_{r-1}-1}\left(\frac{t_{r-1}}{t_r}\right)\log^{k_r-1}\left(\frac{t_r}{z}\right)}{(1-t_1^m)\cdots (1-t_{r-1}^m)(1-t_r^m)}dt_1\cdots dt_r,
\end{align}

Corresponding to Definition \ref{def:1}, we define the multiple zeta function of level $m$ as follows.
\begin{defn} For $k_1,\ldots,k_{r-1}\in\N$ and $\Re(s)>1$, let
\begin{align}\label{d6}
&T^{(m)}_0(k_1,\ldots,k_{r-1},s):= \sum\limits_{1 \le {n_1} <  \cdots  < {n_r}\atop n_i\equiv i\ {\rm mod}\ m} {\frac{1}{{n_1^{{k_1}} \cdots n_{r-1}^{{k_{r-1}}}n_r^s}}}.
\end{align}
Furthermore, as its normalized version, let
\begin{align}\label{d7}
T^{(m)}(k_1,\ldots,k_{r-1},s):=m^r T^{(m)}_0(k_1,\ldots,k_{r-1},s).
\end{align}
\end{defn}
When $k_r>1$, we see that
$${\rm Ath}^{(m)}(k_1,\ldots,k_r;1)=T^{(m)}_0(k_1,\ldots,k_r).$$

\begin{lem} For $k_1,\ldots,k_{r-1}\in\N$ and $\Re(s)>1$,
\begin{align}\label{d8}
T^{(m)}_0(k_1,k_2,\ldots,k_{r-1},s)=\frac{\underbrace{\int\limits_{0}^\infty\cdots\int\limits_{0}^\infty}_{r} \left\{\prod\limits_{j=1}^{r-1}x^{k_j-1}_j\right\}x^{s-1}_r \left\{\prod\limits_{j=1}^{r} \frac{e^{(m-1)(x_j+\cdots+x_r)}}{e^{m(x_j+\cdots+x_r)}-1}\right\}dx_1\cdots dx_r}{\left\{\prod\limits_{j=1}^{r-1}\Gamma(k_j)\right\}\Gamma(s)}.
\end{align}
\end{lem}
\pf The result immediately follows from the definition (\ref{d6}) and the expression $\frac{1}{n^s}=\frac{1}{\Gamma(s)}\int\limits_{0}^\infty e^{-nt}t^{s-1}dt$. \hfill$\square$

\subsection{Some connections between $\psi^{(m)}({\bf k}_r;s)$ and $f_m(z)$}
\begin{defn} For $k_1,\ldots,k_{r}\in\N$ and $\Re(s)>0$, if $f_m(z)$ exist, we can define the level $m$-version of $\xi(k_1,\ldots,k_r;s)$ by
\begin{align}\label{d9}
\psi^{(m)}(k_1,k_2\ldots,k_r;s):=\frac{1}{\Gamma(s)} \int\limits_{0}^\infty \frac{t^{s-1}}{e^t-e^{(1-m)t}}\A^{(m)}\left({k_1,k_2,\ldots,k_r};f_m\left(\frac{t}{m}\right)\right)dt,
\end{align}
where $\A^{(m)}({k_1,k_2,\ldots,k_r};z):=m^r{\rm {Ath}}^{(m)}({k_1,k_2,\ldots,k_r};z)$.
\end{defn}
It is clear that
$$\psi^{(1)}(k_1,k_2\ldots,k_r;s)=\xi(k_1,k_2\ldots,k_r;s)\quad{\rm and}\quad \psi^{(2)}(k_1,k_2\ldots,k_r;s)=\psi(k_1,k_2\ldots,k_r;s).$$

\begin{thm} For $r,k\in\N$, if $f_m(z)$ exist, the following identity holds.
\begin{align}\label{d10}
\psi^{(m)} (\{1\}_{r-1},k;s)&=(-1)^{k-1} \sum\limits_{a_1+\cdots+a_k=r\atop a_1,\ldots,a_k\geq 0} \binom{s+a_k-1}{a_k} T^{(m)}(a_1+1,\ldots,a_{k-1}+1,s+a_k)\nonumber\\
&\quad+\sum\limits_{j=0}^{k-2} (-1)^j T^{(m)}(\{1\}_{r-1},k-j)\cdot T^{(m)}(\{1\}_j,s).
\end{align}
\end{thm}
\pf The method of the proof is similar to that in \cite[Theorem 8]{AM1999} and \cite[Theorem 5.3]{KTA2018}. Given $r,k\geq 1$, introduce the following integrals
\begin{align*}
I^{(r,k)}_{j,m}(s):=\frac{m^{k}}{\Gamma(s)} \underbrace{\int\limits_{0}^\infty\cdots\int\limits_{0}^\infty}_{k-j+1} \frac{\A\left(\{1\}_{r-1},j;f_m\left(\frac{x_j+\cdots+x_k}{m}\right)\right)}{\prod\limits_{l=j}^k \left(e^{x_l+\cdots+x_k}-e^{(1-m)(x_l+\cdots+x_k)} \right)}x^{s-1}_k dx_l \cdots dx_k.
\end{align*}
We compute $I^{(r,k)}_{1,m}(s)$ in two different ways. Firstly, from (\ref{d2}) and (\ref{d4}),
\begin{align*}
\A^{(m)}\left(\{1\}_{r};f_m\left(\frac{x_j+\cdots+x_k}{m}\right)\right)=\frac{(x_1+\cdots+x_r)^r}{r!}.
\end{align*}
Then, by (\ref{d8}), we obtain
\begin{align}\label{d11}
I^{(r,k)}_{1,m}(s)=&\frac{m^{k}}{\Gamma(s)r!} \underbrace{\int\limits_{0}^\infty\cdots\int\limits_{0}^\infty}_{k} \frac{(x_1+\cdots+x_k)^r x^{s-1}_{k}}{\prod\limits_{l=1}^k \left(e^{x_l+\cdots+x_k}-e^{(1-m)(x_l+\cdots+x_k)} \right)}dx_1\cdots dx_k\nonumber\\
=&\frac{m^k}{\Gamma(s)} \sum\limits_{a_1+\cdots+a_k=r \atop a_1,\ldots,a_k\geq 0} \frac{1}{a_1!\cdots a_k!}\underbrace{\int\limits_{0}^\infty\cdots\int\limits_{0}^\infty}_{k}\frac{x_1^{a_1}\cdots x_{k-1}^{a_{k-1}}x_k^{s+a_k-1}}{\prod\limits_{l=1}^k \left(e^{x_l+\cdots+x_k}-e^{(1-m)(x_l+\cdots+x_k)} \right)}dx_1\cdots dx_k\nonumber\\
=& \sum\limits_{a_1+\cdots+a_k=r\atop a_1,\ldots,a_k\geq 0} \binom{s+a_k-1}{a_k} T^{(m)}(a_1+1,\ldots,a_{k-1}+1,s+a_k).
\end{align}

Secondly, by using
\begin{align}\label{d12}
&\frac{\partial}{\partial x_j}{\rm Ath}^{(m)}\left({\{1\}_{r-1},j+1;f_m\left(\frac{x_j+\cdots+x_k}{m}\right)}\right)\nonumber\\
&=m \frac{{\rm Ath}^{(m)}\left({\{1\}_{r-1},j;f_m\left(\frac{x_j+\cdots+x_k}{m}\right)}\right)}{e^{x_j+\cdots+x_k}-e^{(1-m)(x_j+\cdots+x_k)}}
\end{align}
and (\ref{d8}), we compute
\begin{align*}
I^{(r,k)}_{j,m}(s)&=\frac{m^{k}}{\Gamma(s)} \underbrace{\int\limits_{0}^\infty\cdots\int\limits_{0}^\infty}_{k-j+1}\frac{ \frac{\partial}{\partial x_j}\left\{\A\left(\{1\}_{r-1},j+1;f_m\left(\frac{x_j+\cdots+x_k}{m}\right)\right)\right\}}{\prod\limits_{l=j+1}^k \left(e^{x_l+\cdots+x_k}-e^{(1-m)(x_l+\cdots+x_k)} \right)}x^{s-1}_k dx_l \cdots dx_k\nonumber\\
&=T^{(m)}(\{1\}_{r-1},j+1)\cdot T^{(m)}(\{1\}_{k-j-1},s)-I^{(r,k)}_{j+1,m}(s).
\end{align*}
Therefore, using this definition relation repeatedly, we obtain
\begin{align*}
&I^{(r,k)}_{1,m}(s)=\sum\limits_{j=1}^{k-1} (-1)^{j-1} T^{(m)}(\{1\}_{r-1},j+1)\cdot T^{(m)}(\{1\}_{k-j-1},s)+(-1)^{k-1}I^{(r,k)}_{k,m}(s).
\end{align*}
By definition, we have
$$I^{(r,k)}_{k,m}(s)=\psi^{(m)}(\{1\}_{r-1},k;s),$$
and thus
\begin{align}\label{d13}
&I^{(r,k)}_{1,m}(s)=\sum\limits_{j=0}^{k-2} (-1)^{k-j} T^{(m)}(\{1\}_{r-1},k-j)\cdot T^{(m)}(\{1\}_j,s)+(-1)^{k-1}\psi^{(m)}(\{1\}_{r-1},k;s).
\end{align}
Comparing (\ref{d11}) and (\ref{d13}), we obtain the assertion.\hfill$\square$

\begin{thm}\label{thm4.7} For $r,k\in\N$ and $p\in\N_0$, if $f_m(z)$ exist, the following identity holds.
\begin{align}\label{d14}
\psi^{(m)}(\{1\}_{r-1},k;p+1)=\sum\limits_{a_1+\cdots+a_k=p\atop a_,\ldots,a_k\geq 0} \binom{a_k+r}{r} T^{(m)}(a_1+1,\ldots,a_{k-1}+1,a_k+r+1).
\end{align}
\end{thm}
\pf By (\ref{d12}), we have
\begin{align*}
&\psi^{(m)}(\{1\}_{r-1},k;p+1)\nonumber\\&=\frac{m}{\Gamma(p+1)} \int\limits_{0}^\infty \frac{t^{p}}{e^t-e^{(1-m)t}}\A^{(m)}\left(\{1\}_{r-1},k;f_m\left(\frac{t}{m}\right)\right)dt\nonumber\\
&=\frac{m^2}{\Gamma(p+1)} \int\limits_{0}^\infty \frac{t^{p}_k}{e^{t_k}-e^{(1-m)t_k}}\int\limits_{0}^{t_k}\frac{\A^{(m)}\left(\{1\}_{r-1},k-1;f_m\left(\frac{t_{k-1}}{m}\right)\right)}{e^{t_{k-1}}-e^{(1-m)t_{k-1}}}dt_{k-1}dt_k\nonumber\\
&=\frac{m^k}{p!} \underbrace{\int\limits_{0}^\infty \int\limits_{0}^{t_k}\cdots\int\limits_{0}^{t_2}}_{k} \frac{t^p_k {\A}^{(m)}\left(\{1\}_r;f_m\left(\frac{t_1}{m}\right)\right)}{(e^{t_k}-e^{(1-m)t_k})\cdots (e^{t_1}-e^{(1-m)t_1})}dt_1dt_2\cdots dt_k\\
&=\frac{m^k}{p!r!} \underbrace{\int\limits_{0}^\infty \int\limits_{0}^{t_k}\cdots\int\limits_{0}^{t_2}}_{k} \frac{t^p_k t^r_1}{(e^{t_k}-e^{(1-m)t_k})\cdots (e^{t_1}-e^{(1-m)t_1})}dt_1dt_2\cdots dt_k.
\end{align*}
By the change of variables
$$t_1=x_k,t_2=x_{k-1}+x_k,\ldots,t_k=x_1+\cdots+x_k,$$
we obtain
\begin{align*}
\psi^{(m)}(\{1\}_{r-1},k;p+1)=&\frac{m^k}{p!r!} \underbrace{\int\limits_{0}^\infty \int\limits_{0}^{t_k}\cdots\int\limits_{0}^{t_2}}_{k} \frac{(x_1+\cdots+x_k)^p x_k^r}{\prod\limits_{l=1}^k (e^{x_l+\cdots+x_k}-e^{(1-m)(x_l+\cdots+x_k)})}dt_1dt_2\cdots dt_k\\
=&\sum\limits_{a_1+\cdots+a_k=p\atop a_,\ldots,a_k\geq 0} \binom{a_k+r}{r} T^{(m)}(a_1+1,\ldots,a_{k-1}+1,a_k+r+1).
\end{align*}
Thus, the proof of Theorem \ref{thm4.7} is finished.\hfill$\square$

\begin{cor} For $r,k\in\N$, if $f_m(z)$ exist, then we have the ``height one" duality
\begin{align}\label{d15}
T^{(m)}(\{1\}_{r-1},k+1)=T^{(m)}(\{1\}_{k-1},r+1).
\end{align}
\end{cor}
\pf Setting $p=0$ in (\ref{d14}) gives
\begin{align}\label{d16}
\psi^{(m)}(\{1\}_{r-1},k;1) =T^{(m)}(\{1\}_{k-1},r+1).
\end{align}
In general, from the definition, we have
\begin{align*}
\psi^{(m)}(k_1,\ldots,k_r;1)&=m \int\limits_{0}^\infty \frac{{\A}^{(m)}\left(k_1,\ldots,k_r;f_m(t/m)\right)}{e^t-e^{(1-m)t}}dt\\
&=T^{(m)}(k_1,\ldots,k_{r-1},k_r)
\end{align*}
and in particular
\begin{align}\label{d17}
\psi^{(m)}(\{1\}_{r-1},k;1) =T^{(m)}(\{1\}_{r-1},k+1).
\end{align}
Thus, from (\ref{d16}) and (\ref{d17}) we obtain (\ref{d15}).\hfill$\square$

\subsection{Duality relation for $\psi^{(m)}({\bf k}_r;p+1)$}

If $f_m(z)$ exist, from (\ref{d2}), we have
$$\A^{(m)}(\{1\}_r;f_m(t/m))=\frac{t^r}{r!}.$$
Hence,
\begin{align*}
\psi^{(m)}(k_1,\ldots,k_r;p+1)&=m \int\limits_{0}^\infty \frac{\A^{(m)}\left(\{1\}_p;f_m(t/m)\right)\A^{(m)}\left(k_1,\ldots,k_r;f_m(t/m)\right)}{e^t-e^{(1-m)t}}dt.
\end{align*}
By the change of variable $u=f_m(t/m)$, we obtain
\begin{align}\label{d18}
\psi^{(m)}(k_1,\ldots,k_r;p+1)&= \int\limits_{0}^1 \frac{\A^{(m)}(\{1\}_p;u)\A^{(m)}(k_1,\ldots,k_r;u)}{u}du.
\end{align}
\begin{thm}\label{thm4.9} For positive integers $p,q,r\in\N$ and $k_1,\ldots,k_r\in \N\setminus\{1\}$, if $f_m(z)$ exist, then the following duality relation holds:
\begin{align}\label{d19}
&\psi^{(m)}(\{1\}_{q-1},k_1,\ldots,k_{r-1},k_r-1;p+1)-(-1)^{k_1+k_2+\cdots+k_r}\psi^{(m)}(\{1\}_{p-1},k_r,\ldots,k_{2},k_1-1;q+1)\nonumber\\
&=\sum\limits_{j=0}^{r-1} (-1)^{\mid\stackrel{\leftarrow}{{\bf k}}_j\mid}\sum\limits_{i=1}^{k_{r-j}-2} (-1)^{i-1} T^{(m)}(\{1\}_{p-1},k_r,k_{r-1},\ldots,k_{r+1-j},i+1)\nonumber\\&\quad\quad\quad\quad \quad\quad \quad\quad \quad\quad \times T^{(m)}(\{1\}_{q-1},k_1,k_{2},\ldots,k_{r-j-1},k_{r-j}-i)\nonumber\\
&\quad+\sum\limits_{j=0}^{r-2} (-1)^{\mid\stackrel{\leftarrow}{{\bf k}}_{j+1}\mid} \lim_{x\rightarrow 1}\left\{\begin{array}{l} \A^{(m)}(\{1\}_{p-1},k_r,k_{r-1},\ldots,k_{r-j};x)\\ \quad\times\A^{(m)}(\{1\}_{q-1},k_1,k_2,\ldots,k_{r-j-1},1;x)\\ -\A^{(m)}(\{1\}_{p-1},k_r,k_{r-1},\ldots,k_{r-j},1;x)\\ \quad\times\A^{(m)}(\{1\}_{q-1},k_1,k_2,\ldots,k_{r-j-1};x)\end{array}\right\}.
\end{align}
\end{thm}
\pf By considering
\begin{align}\label{d20}
&\psi^{(m)}(\{1\}_{q-1},k_1,\ldots,k_{r-1},k_r-1;p+1)\nonumber
\\&=\int\limits_{0}^1 \frac{\A^{(m)}(\{1\}_p;u)\A^{(m)}(\{1\}_{q-1},k_1,\ldots,k_{r-1},k_r-1;u)}{u}du.
\end{align}
Then, by applying the same arguments as in the proof of Theorem \ref{thm3.3}, we may easily deduce the result.\hfill$\square$

\begin{re} It is possible to obtain the explicit evaluation of the limit in (\ref{d19}) by using a similar method as in the proof of Theorem \ref{thm3.3}.

\end{re}
{\bf Acknowledgments.} The authors express deepest gratitude to their supervisor Professor Masanobu Kaneko for
his valuable comments and encouragement. The corresponding author Ce Xu is supported by the National Natural Science Foundation of China (Grant No. 12101008), the Natural Science Foundation of Anhui Province (Grant No. 2108085QA01) and the University Natural Science Research Project of Anhui Province (Grant No. KJ2020A0057).

\medskip

\noindent{\bf Disclosure statement.} The authors report there are no competing interests to declare.

 \end{document}